\documentclass{article}

\usepackage{PRIMEarxiv}

\usepackage[utf8]{inputenc} % allow utf-8 input
\usepackage[T1]{fontenc}    % use 8-bit T1 fonts
\usepackage{hyperref}       % hyperlinks
\usepackage{url}            % simple URL typesetting
\usepackage{booktabs}       % professional-quality tables
\usepackage{amsfonts}       % blackboard math symbols
\usepackage{nicefrac}       % compact symbols for 1/2, etc.
\usepackage{microtype}      % microtypography
\usepackage{lipsum}
\usepackage{fancyhdr}       % header
\usepackage{graphicx}       % graphics
\graphicspath{{media/}}     % organize your images and other figures under media/ folder

 \usepackage{mathtools}
 \usepackage{cite}
 \usepackage{amsmath}
\usepackage{amssymb}
\usepackage{bm}

 \usepackage{algorithm}
 \usepackage{algpseudocode}
 \usepackage{eqlist}
\usepackage[overload]{empheq}
\usepackage{mathtools}
\usepackage{xcolor}

 \newtheorem{thm}{Theorem}

\newtheorem{rem}{Remark}
\newtheorem{lem}{Lemma}

\newcommand*{\QEDA}{\null\nobreak\hfill\ensuremath{\blacksquare}}%

\title{Stability of Finite Receding Horizon Control: A Complementary Approach
\thanks{Corresponding author: Wen-Hua Chen.} 
}

\author{
Wen-Hua Chen \\
 Department of Aeronautical and Automotive Engineering  \\
 Loughborough University  \\
 Loughborough LE11 3TU, U.K.\\
  \texttt{w.chen@lboro.ac.uk} \\
   \And
 Yunda Yan\\
 School of Engineering and Sustainable Development \\
  De Montfort University  \\
  Leicester, LE1 9BH, U.K.\\
  \texttt{yunda.yan@dmu.ac.uk}  \\
}

\begin{document}
\maketitle

\begin{abstract}
This paper presents a complementary approach to establish stability of finite receding horizon control 	with a terminal cost. First a new augmented stage cost is defined by rotating the terminal cost. Then a one-step optimisation problem is defined based on this augmented stage cost. It is shown that a slightly modified Model Predictive Control (MPC) algorithm is stable if the value function of the augmented one-step cost (OSVF) is a control Lyapunov function. The proposed stability condition is completely complementary to the existing terminal cost based MPC stability conditions in the sense that they are mutually excluded with each other. By using this approach, we are able to establish stability for MPC algorithms with zero terminal cost or even negative terminal cost as special cases. Combining this new approach with the existing MPC stability theory, we are able to significantly relax the stability requirement on MPC and extend the design space where stability are guaranteed. The proposed approach will help to further reduce the gap between stability theory and practical applications of MPC and other optimisation based control methods.  
\end{abstract}

% keywords can be removed

\section{Introduction}
 It is well known that stability of Receding Horizon Control (RHC), also known as Model Predictive Control (MPC), with a finite horizon can not be guaranteed. The most widely used approach to guarantee the stability of a receding horizon controller is to add a terminal cost, $V_f(x_{k+N|k})$, in its cost function where $x_{k+N|k}$ is the terminal state in the prediction horizon \cite{mayne2000constrained,CheBalOre00,fontes2001general,primbs1999nonlinear}. Consequently, a well-established sufficient condition is that the system $x^+=f(x,u)$ is stable under an MPC algorithm if there exists \emph{any} control $u$ in the admissible set $\mathbb{U}$ such that 
\begin{equation} \label{eq:original}
	-V_f(x)+ l (x,u)+V_f(f(x,u))<0
\end{equation} 
 for any $x \neq 0$ where $l(x,u)$ is the stage cost.  That is, a terminal cost $V_f(x_{k+N|k})$ that covers the time-go-cost (i.e. from the terminal state to the final equilibrium) shall be added into the cost function for on-line optimisation. A natural question arising is, what happens if condition (\ref{eq:original}) is violated, or more precisely, if there does not exist an admissible control $u \in \mathbb{U}$ such that condition (\ref{eq:original}) holds. Formally, this can be  represented as, for \emph{all} $u \in \mathbb{U}$ and {$x \neq 0$},     
\begin{equation} \label{eq:newdefl}
	\ell(x,u,x^+) \coloneqq -V_f(x)+ l (x,u)+V_f(f(x,u))>0
\end{equation} 
This question is more interesting since it includes the case that $V_f(x)=0$, i.e. no terminal cost. In the case of zero terminal cost, condition (\ref{eq:newdefl}) is satisfied as long as the stage cost, $l(x,u)$, is positive definite. Unfortunately, numerous counterexamples exist where condition (\ref{eq:newdefl}) holds but the resulted MPC is unstable, {e.g. \cite{kouvaritakis2002letter, chen2020model}}. It is noted that significant effort and progress have been made in relaxing the requirement on the terminal cost in satisfying (\ref{eq:original}) \cite{jadbabaie2005stability,grimm2005model,GruPan17,faulwasser2018economic}. For example, stability can still be achieved if the receding horizon is sufficiently long even in the case when the terminal weight is zero \cite{grimm2005model}. So far, there is little work devoted to directly  exploring the space defined by condition (\ref{eq:newdefl}) which is \emph{opposite} to the design space defined by (\ref{eq:original}).  

The motivation of this paper is to answer the question \emph{under what condition a finite receding horizon controller satisfying (\ref{eq:newdefl}) is stable}. Condition (\ref{eq:newdefl}) is equivalent to  $m(x)>0$, where 
\begin{equation} \label{eq:m(x)}
    m(x)\coloneqq \min_{u\in \mathbb{U}} ~\ell(x,u,x^+)
\end{equation} 
is referred as the One-Step Value Function (OSVF). That plays a key role in establishing stability of MPC in this paper. The main contribution of this paper is to present a slightly modified MPC algorithm which is stable \emph{as long as the corresponding OSVF is a control Lyapunov function (CLF)}.

Since condition (\ref{eq:newdefl}) is precisely opposite to condition (\ref{eq:original}) used in establishing stability in the current MPC literature for more than two decades \cite{mayne2000constrained,lee2011model,GruPan17},  our results provide a \emph{complementary} approach to the existing terminal cost based MPC framework. That is, all cases covered by our new conditions do not satisfy the existing stability conditions developed based on (\ref{eq:original}). Furthermore, our conditions allow stability to be established if there is no terminal cost or even a negative terminal cost (see numerical examples). It clearly shows that, despite the huge success in developing stability theory for MPC in the last three decades \cite{mayne2000constrained,rawlings2017model,faulwasser2018economic}, there is still a significant, unexplored, design space where stability is possible to achieve by optimising a finite horizon cost function. By combining the new analysis approach in this paper with the existing ones, particularly the existing terminal cost based MPC stability conditions, we are in a much stronger position in analysing stability of MPC and other optimisation based control methods and designing stability guaranteed optimisation based algorithms with a much less restricted freedom and space. This is the main motivation behind this paper.  

There are two key ideas behind the proposed approach. The first is to define a new stage cost, $\ell(x,u,x^+)$, as in (\ref{eq:newdefl}) by rotating the terminal cost. 
%This can be considered as an extension of the augmented stage cost defined in \cite{chen2020model} where the traditional stage cost $l(x,u)$ is replaced by $\ell(x^+,u)$. Stability of MPC without terminal cost is established using this new defined stage cost. These two stage costs, i.e. $l(x,u)$ and $\ell(x^+,u)$, could be considered as two special cases of the augmented stage cost $\ell(x,u,x^+)$ in this paper although there is a subtle difference. Furthermore, 
{It shall be pointed out that the stage cost in (\ref{eq:newdefl}) is somehow similar to the rotation stage cost technique widely used in economic MPC (EMPC) \cite{faulwasser2018economic}. To establish stability of EMPC, the stage cost is first rotated by a storage function, $\lambda(x)$, (normally being non-negative \cite{faulwasser2018economic}). Subsequently, it is required that the modified terminal cost, $\tilde{V_f}(x)=V_f(x)+\lambda(x)$,  is a CLF. This is the same technique as in the current tracking MPC  (\ref{eq:original}); for example, see \cite{faulwasser2018economic,angeli2011average,diehl2010lyapunov}. It is interesting to notice that the widely used dissipativity inequality is equivalent to condition (\ref{eq:newdefl}) when the storage function $\lambda(x)$ is specifically chosen as $-V_f(x)$. However, we cannot derive any stability results under this specific choice as discussed in Remark \ref{remEMPC}.}

The second idea is that it is observed that the change of the initial state related cost in a cost function involved in online optimisation involved in MPC does not change the optimal control sequence and the MPC algorithm.
%\cite{chen2020model}.
By combining these two ideas, we are able to modify the MPC cost function but without changing its behaviour. \emph{Fundamentally, when the optimal value function of the modified cost is used as a Lyapunov function candidate to establish its stability, the combination of these two ideas enables a wide range of possible Lyapunov function candidates to be used to establish stability for the same MPC algorithm. This breaks down the rigid link between the cost function used in optimisation in an MPC algorithm and the optimal value function used as a Lyapunov function candidate in the current approaches} \cite{mayne2000constrained,lee2011model}. { The fact that changing initial state cost does not alter the optimisation involved in MPC has been exploited in \cite{faulwasser2015design,GruPan17}. By further exploiting this property in this paper, we are able to explore new design space for MPC stability that was not possible previously.}   

Technically, to establish stability by the virtue of the OSVF being a CLF, we have to slightly modify the MPC algorithm. That is, rather than using a fixed terminal constraint as in the current terminal cost based  framework, we construct a contractive terminal set by making use of the  OSVF being a CLF.   The idea of using a contractive terminal set can be found in \cite{de2000contractive}; however, the whole optimal control sequence (not only the first element) is applied into the system in \cite{de2000contractive}. This does not conform with the MPC setup in this paper.

This paper is organised as follows. In Section \ref{sec2}, a new stage cost is proposed by rotating the  terminal cost in the conventional cost function without affecting the optimal trajectories. Based on this new augmented stage cost and some mild assumptions, a new MPC algorithm with stability guarantee is proposed in  Section \ref{sec3}. In Section \ref{sec4}, a procedure for finding a terminal cost satisfying the new stability conditions is developed.  Results are provided by two illustrative examples, giving further insight into the proposed approach, in Section  \ref{sec5}. Finally, Section \ref{sec6} concludes the paper.

Notation. $\mathbb{I}$  and $\mathbb{R}$ are integers and real numbers, respectively, where subscripts may be added to give specific ranges.  For any vector $x\in \mathbb{R}^n$, {$|x|$ denotes the 2-norm} and $|x|_P$ is defined by $|x|_P\coloneqq x^TPx$, where $P\in \mathbb{R}^{n\times n}$ is a symmetric matrix, but not necessary to be positive-definite here. The superscript $^*$ is used to indicate the optimal solution,  the optimal control sequence or state under the optimal control sequence.

\section{Augmented stage cost and equivalent MPC problems}\label{sec2}
 
Consider a nonlinear discrete-time system described by
\begin{equation}\label{sys}
	\begin{aligned}
x_{k+1}=f(x_k,u_k),~k\in \mathbb{I}_{\ge0}
	\end{aligned}
\end{equation} 
 or, more succinctly $x^+=f(x,u)$,
where $x \text { or }x_k\in \mathbb{X}$ is the current state, $u\text { or }u_k\in \mathbb{U}$ is the current control, and $x^+\text { or }x_{k+1}$ is the successor state. 
It  is assumed in this paper that
\begin{description}
	\item[{(A1)}] {$f(x,u): \mathbb{X}\times \mathbb{U} \to \mathbb{X}$ is  continuous  and $f(0,0)=0$. $\mathbb{X}\subseteq \mathbb{R}^n$ and $\mathbb{U}\subseteq \mathbb{R}^m$ are closed and  contain the origin in the interior. }  
\end{description}

We start from a classic finite horizon optimisation problem, with the cost function
 \begin{equation}\label{J1}
J (x_k, \bm{\mathrm u}_k) \coloneqq\sum_{i=0}^{N-1} l(x_{k+i|k},u_{k+i|k})+V_f(x_{k+N|k}),
\end{equation}
where $x_{k|k}=x_k$ and $x_{k+i+1|k}=f(x_{k+i|k},u_{k+i|k}),~i\in \mathbb{I}_{0:N-1}$; the control and state sequences are\protect\footnotemark[1] $\bm{\mathrm u}_k\coloneqq \left(u_{k|k}, \cdots, u_{k+N-1|k}\right)$ and $\bm{\mathrm x}_k\coloneqq \left(x_{k|k},\cdots, x_{k+N|k}\right)$, subject to the constraints $\mathbb{U}$ and $\mathbb{X}$, $N\ge 1$ is the length of the control horizon, $l(x,u)$ and  $V_f(x)$ are  the stage cost and the terminal cost respectively. This satisfies the following assumption.
\begin{description}
	\item[(A2)] $l(x,u): \mathbb{X}\times \mathbb{U} \to \mathbb{R}$ and $V_f(x): \mathbb{X} \to \mathbb{R}$ are continuous. Furthermore, $l(0,0)=0$ and $V_f(0)=0$.
\end{description}

\footnotetext[1]{For simplicity,  when used in algebraic  expressions, $\bm{\mathrm x}_k$ denotes the column vector $\left(x^T_{k|k}, \cdots, x^T_{k+N|k}\right)^T$; similarly  in algebraic expressions of $\bm{\mathrm u}_k$. }

It shall be pointed out that there is no requirement that $(l,x)$ or $V_f(x)$ shall be positive definite in this paper.

Following the conventional MPC framework, the first control action of the optimal control sequence $\bm{\mathrm u}_k^*$ is applied into the system, $u_k=u^*_{k|k}$,
which implies that
\begin{equation}\label{ns}
	x_{k+1}=f(x_k,u_k)=f\left(x_{k|k},u^*_{k|k}\right)=x^*_{k+1|k}.
\end{equation}
For the widely used cost function in (\ref{J1}), it is observed that the cost  that is \textit{only} associated with the initial state, $x_{k|k}$, is never changed by any control sequence. That is,  this term is completely independent from the optimisation process.

Motivated by this intuitive but important observation, we rewrite the cost function as 
\begin{equation}\label{lp}
\begin{aligned}
&J (x_k, \bm{\mathrm u}_k)\\
&= \sum_{i=0}^{N-1} l(x_{k+i|k},u_{k+i|k})+V_f(x_{k+N|k})\\
&=V_f(x_{k|k})\underbrace{-V_f(x_{k|k})+l(x_{k|k},u_{k|k})+V_f(x_{k+1|k})}_{\ell(x_{k|k},u_{k|k},x_{k+1|k})}\\ 
&~~~~~~\vdots\\ 
&~~~\underbrace{-V_f(x_{k+N-1|k})+l(x_{k+N-1|k},u_{k+N-1|k})  +V_f(x_{k+N|k}) }_{\ell(x_{k+N-1|k},u_{k+N-1|k},x_{k+N|k})}\\
&=V_f(x_{k|k})+\sum_{i=0}^{N-1} \ell(x_{k+i|k},u_{k+i|k},x_{k+i+1|k}), 
\end{aligned}
\end{equation}
where $\ell(x,u,x^+)$  is the new augmented stage cost  and has been defined in (\ref{eq:newdefl}) by rotating the terminal cost. It shall be highlighted that the information of the system dynamics is embedded in this new stage cost.  By ignoring the initial one, $V_f(x_{k|k})$, since it does not affect the solution of the optimisation, we can define the new cost function for MPC online optimisation as
\begin{equation}\label{J2}
	\mathcal{J} (x_k, \bm{\mathrm u}_k)\coloneqq  \sum_{i=0}^{N-1}\ell(x_{k+i|k},u_{k+i|k},x_{k+i+1|k}).
\end{equation}
Considering the optimisation problems based on different cost functions, $J(x_k, \bm{\mathrm u}_k)$ and $\mathcal{J} (x_k, \bm{\mathrm u}_k)$, with the same initial state $x_k$ and the same constraints, the optimal trajectories should be the same, regardless of the different value functions. That is
\begin{equation}\label{JJ}
	\begin{aligned}
\min_{\bm{\mathrm u}_k} \mathcal{J} (x_k, \bm{\mathrm u}_k)&= -V_f(x_k)+\min_{\bm{\mathrm u}_k} J (x_k, \bm{\mathrm u}_k)\\
	\arg \min_{\bm{\mathrm u}_k} \mathcal{J} (x_k, \bm{\mathrm u}_k)&=  \arg \min_{\bm{\mathrm u}_k} J (x_k, \bm{\mathrm u}_k).
	\end{aligned}
\end{equation}
Considering (\ref{JJ}), the design and analysis on the conventional MPC  are equivalent to that with a new cost function $\mathcal{J} (x_k, \bm{\mathrm u}_k)$, which is given in a form without any terminal cost.

\begin{rem} \label{remEMPC}
{The stage cost rotation technique has been widely used for EMPC in establishing its stability \cite{angeli2011average,diehl2010lyapunov,zanon2014indefinite,zanon2018economic}. More specifically, a dissipativity inequality is defined with a non-negative storage function $\lambda(x)$ (see Definition 3.1 in \cite{faulwasser2018economic}), that is, \begin{equation}
    -\lambda(x)+\lambda(f(x,u)) -l(x,u) \le -l(x_s,u_s)
\end{equation}
is satisfied for all $x$ and $u$, where $x_s$ and $u_s$ are in steady state operating condition and $l(x_s,u_s)$ can be treated as zero when the equilibrium point is at the origin. 
Then stability of EMPC can be established by mainly resorting to standard stabilising MPC theory. That is, to establish stability, the modified terminal cost $\tilde{V}_f(x)=V_f(x)+\lambda(x)$ is required to be a Control Lyapunov Function (CLF), satisfying condition (\ref{eq:original}). The authors are aware that  significant effort and progress has been made in further developing stability analysis of EMPC with various extensions (see, for example, \cite{dong2017generalized,zanon2014indefinite,zanon2018economic,faulwasser2018asymptotic,faulwasser2018economic}).

Condition (\ref{eq:newdefl}) can be interpreted as a dissipativity inequality with the corresponding storage function as $\lambda(x)=-V_f(x)$. It is interesting to notice that the resulted storage function $\lambda(x)$ is quite often negative definite since the terminal cost is, in general, positive definite, which makes it different from a standard dissipation definition. Most importantly, with the choice of $\lambda(x)=-V_f(x)$, the modified terminal cost becomes $\tilde{V}_f(x)=0$, which does not satisfy the requirement that the terminal cost must be a CLF, i.e. condition (\ref{eq:original}). We are not able to assert stability of EMPC or any other MPC schemes by following this approach. 

Furthermore, based on the dissipativity condition, similar functions have recently been defined in EMPC using a general storage function which is referred to as a control dissipativity function in \cite{lazar2021stabilization} and control storage function in \cite{dong2017generalized}. }

\end{rem}

{\begin{rem}
The approach proposed in this paper is  applicable to both tracking MPC and EMPC directly. It follows from Assumption \textbf{(A3)} that the stage cost $l(x,u)$ could be positive or negative as long as condition (\ref{eq:newdefl}) is satisfied. In that sense, we also present an alternative stability condition for EMPC where the rotated terminal cost $\tilde{V_f}(x)$ does not satisfy the requirements in the current stability conditions, e.g. to be a CLF.
\end{rem}}

 \section{MPC algorithm and stability}\label{sec3}

\subsection{New terminal weight based MPC algorithm}

Here we define an auxiliary optimisation problem to minimise the augmented stage cost $\ell(x,u,x^+)$ in (\ref{eq:newdefl}) or (\ref{lp}) 
\begin{equation}\label{m}
\begin{aligned}
m(x)\coloneqq &\min_{u\in \mathbb{U}} ~\ell(x,u,x^+)\\
	=&\min_{u\in \mathbb{U}}~-V_f(x)+l(x,u)+V_f(x^+)\\ 
	&\text{s.t.}~x\in \Omega \subseteq  \mathbb{X},~u\in \mathbb{U},~f(x,u) \in \Omega,
\end{aligned}
\end{equation}
where $m(x): \Omega \subseteq \mathbb{X}\mapsto \mathbb{R}$ is
 the One-Step  Value Function (OSVF)  of the augmented stage cost and $\Omega$ is a control  invariant set.
 
 It shall be highlighted that $m(x)>0$ iff condition (\ref{eq:newdefl}) holds. In other words, $m(x)>0$ only if there does not exist \emph{any} control $u \in \mathbb{U}$ such that condition (\ref{eq:original}) holds. We will investigate the design space where $m(x)>0$ is positive definite. This is entirely complementary to the design space investigated by the existing MPC stability theories and has not been explored previously. 
 
We  define the sublevel set of $m(x)>0$ as
	\begin{equation}
		\Omega(\alpha)\coloneqq \left\{ x \in \mathbb{R}^n: m(x) \leq \alpha\right\}\subseteq  \Omega \subseteq  \mathbb{X},
	\end{equation}
where  $\alpha$ is a positive constant. As  briefly mentioned  in the introduction, the proposed MPC in this paper  is based on the condition that OSVF is a CLF. That is, the following assumption is imposed in this paper.
\begin{description}
	\item[(A3)] 
The OSVF $m(x)$ is a  CLF  for system (\ref{sys}) with respect to a set $\Omega({\alpha}_0)$, i.e., for any $x\in\Omega({\alpha}_0)$,
there exist a control $u\in \mathbb{U}$ and {$\mathcal{K}_\infty$ functions $\beta_1(\cdot),\beta_2(\cdot),\beta_3(\cdot)$}  such that
\begin{align}
{\beta_1 (|x|)	\leq m(x)\leq \beta_2 (|x|)\label{mv1}}\\
{m(f(x,u)) - m(x)\leq -\beta_3 (|x|),\label{mv2}}
\end{align}
where  $\alpha_0>0$.
\end{description}

Based on Assumption \textbf{(A3)}, we are now able to construct
a contractive terminal set, which is the only difference from the conventional MPC with terminal cost. The new algorithm is  described below.

\textbf{Algorithm: }
\begin{quote}
\begin{eqlist*}[\eqliststarinit\def\makelabel#1{#1}\labelsep1em]
 \item[\textbf{Offline:}\\\textbf{Step 1}] For the given system  (\ref{sys}) subject to  the constraints $\mathbb{X}$ and $\mathbb{U}$,  define the stage cost $l(x,u)$ and the terminal cost $V_f(x)$, and calculate  terminal constraint $\alpha_0$. Initialise the time as $k=0$.
       \item[\textbf{Online:}\\\textbf{Step 2}] At time $k$, measure the current state $x_k$. Solve the following optimisation problem 
         \begin{equation}\label{op} 
     \begin{aligned}
     	&V(x_k, {\alpha_k})\coloneqq		\min_{\bm{\mathrm u}_k} \mathcal{J} (x_k, \bm{\mathrm u}_k)\\
     	&=\min_{\bm{\mathrm u}_k}\sum_{i=0}^{N-1}\ell\left(x_{k+i|k},  u_{k+i|k},x_{k+i+1|k}\right)\\
     		&\text{s.t. } x_{k|k}=x_k\\
     		& ~~~~~ x_{k+i+1|k}=f(x_{k+i|k},u_{k+i|k})\\
     		&~~~~~ x_{k+i|k}\in \mathbb{X}, u_{k+i|k}\in \mathbb{U},~i\in   \mathbb{I}_{0:N-1}\\
     		&~~~~~x_{k+N|k}\in   \Omega(\alpha_k)  
     \end{aligned}
\end{equation}
and obtain  the optimal solution $\bm{\mathrm u}^*_k$ and $\bm{\mathrm x}^*_k$. 
 Apply the optimal control $u_k=u^*_{k|k}$ to  system (\ref{sys}).
 
   \item[\textbf{Step 3}]    Denote the minimum between  $m\left(x^*_{k+N|k}\right)$ and $m\left(x^*_{k+1|k}\right)$ as
   \begin{equation}\label{mm}
   \begin{aligned}
   		 m \left(x^*_{k+s|k}\right) \coloneqq\min\left\{ m\left(x^*_{k+1|k}\right),  m\left(x^*_{k+N|k}\right)\right\}, 	s\in\{1,N\}.
   \end{aligned}
   \end{equation}
{Update the terminal set by the following rule
\begin{equation}\label{am}
	\alpha_{k+1} =	\left\{
	\begin{aligned}
	    &m \left(x^*_{k+s|k}\right)-\delta,  \text{ if } m \left(x^*_{k+s|k}\right)\ge\delta\\ 
	    &0, \text{ if } m \left(x^*_{k+s|k}\right)<\delta,
	\end{aligned}\right.
\end{equation}
where $\delta>0$ is a sufficient small constant. }
    \item[\textbf{Step 4}] $k\leftarrow k+1$ and go to \textbf{Step 2}.
\end{eqlist*}
\end{quote}

%p&\in \left(0,\frac{\beta_1}{\beta_3}\right)\cap(0,1]

%  \emph{it may make it more rigorous mathematically but brings more problem. We have to calculate $\beta_1$ and $\beta_3$ which may be not easy before starting the algorithm, this assumption shall be implied by A4!\\
%  Yunda: I avoid it by introducing a sufficient small constant $\gamma$.}

\begin{rem}
Compared with the conventional MPC setup,  Step 3 is the \textit{only} difference. To establish stability, there are also three key terminal elements in this algorithm: terminal cost, terminal control and terminal constraints. Instead of using a fixed terminal constraint, a contractive terminal constraint is constructed in Step 3 where one-step ahead optimisation is involved in (\ref{mm}) and (\ref{am}). It will be shown that under the condition that, for a given terminal cost, if the corresponding OSVF is a CLF, then it is always feasible to construct a contractive terminal constraint. That is, Step 3 and the new MPC algorithm is always feasible as long as it is feasible at the beginning. We will discuss how to choose a terminal cost such that the corresponding OSVF is a CLF in Section~\ref{sec4}. 
%On the other hand,  we only need to choose a sufficiently small gain $\gamma$ without considering any coefficients associated with the CLF condition in Assumption \textbf{(A4)}.  
%These two facts make the proposed MPC algorithm very friendly for practical application.
\end{rem}

\begin{rem}
One of the drawbacks in the current terminal weight based MPC framework is how to strike a good balance between performance and stability, particularly for nonlinear or uncertain system.  For a given stage cost, it is quite difficult to estimate the optimal cost-to-go for these systems. Since, in order to achieve stability, it requires the terminal cost to cover the cost to go,  quite often a conservative terminal cost function has to be employed to ensure the unknown optimal value function is covered. This may lead to  poor performance. The total cost to be optimised in MPC consists of the summed stage cost and the terminal cost. If the chosen terminal cost is much larger than the summed stage cost, the optimisation process actually puts more effort on the terminal cost since the former carries a higher relative weight. This is despite the fact that the stage cost may represent the performance requirements better. Our approach avoids this problem by exploring the design space where the terminal cost cannot cover the cost-to-go, even if it is zero. Our hypothesis is that this may give a better MPC performance for nonlinear or uncertain systems, as shown in the second example in Section~\ref{sec5}.          
\end{rem}
The main result is given by the following theorem while the detailed proof is given by the subsequent sections.
{ 
\begin{thm}\label{t1}
Suppose that 
\begin{itemize} 
  \item  Assumptions \textbf{(A1)}-\textbf{(A3)} are satisfied,
  \item the optimisation problem (\ref{op}) is feasible at time $k=0$,
  \item the parameter $\alpha_0$ is chosen such that the terminal set $\Omega(\alpha_0) \subseteq \Omega$
\end{itemize}
Then, the closed-loop system with the proposed MPC is recursively feasible and,  eventually, asymptotically stable.
%after time $\alpha_0/\delta$.
\end{thm}
}
%p&\in \left(0,\frac{\beta_1}{\beta_3}\right)\cap(0,1]

\subsection{Recursive  feasibility} \label{sec:RF}
Recursive  feasibility of the proposed MPC is mainly based on the condition of OSVF being a CLF, seeing Assumption \textbf{(A3)}.

{
\textit{Proof of  recursive  feasibility in Theorem \ref{t1}.} Supposing that the optimisation (\ref{op}) is feasible at time $k \in  \mathbb{I}_{\ge0}$,   the optimal control and state sequences exist, which can be denoted as 
\begin{align*}
	\bm{\mathrm u}^*_k&=\left(u^*_{k|k}, \cdots, u^*_{k+N-1|k}\right)\\
	\bm{\mathrm x}^*_k&=\left(x^*_{k|k}, \cdots, x^*_{k+N|k}\right),
\end{align*}
where the terminal state satisfies that 
\begin{equation}\label{tc}
	 x^*_{k+N|k}\in   \Omega(\alpha_k) \Leftrightarrow m\left(x^*_{k+N|k}\right)\leq \alpha_k .
\end{equation}
In what follows, we will construct the feasible control and state sequences at time $k+1$ based on the current optimal sequences. The case when the terminal set shrinks to the origin is ignored as it is degenerated into the conventional MPC with an equality terminal constraint. It is straightforward to show its recursive feasibility. To make the discussion clear, we use the superscripts, $a$ and $b$, to represent two cases due to the different results of (\ref{mm}).
}

{
In the case of $m\left(x^*_{k+1|k}\right)\ge m\left(x^*_{k+N|k}\right)$, based on (\ref{am}), if we choose $\delta\in  (0,\beta_3(|x^*_{k+N|k}|)]$,  then
\begin{equation}\label{a2}
\begin{aligned}
		\alpha_{k+1}&=m\left(x^*_{k+N|k}\right)-\delta\\
		&\geq m\left(x^*_{k+N|k}\right)-\beta_3\left(\left|x^*_{k+N|k}\right|\right).
\end{aligned}
\end{equation}
By Assumption \textbf{(A3)} and since  $x^*_{k+N|k}\in   \Omega(\alpha_k)$,  there  exists a control $u^a_{k+N|k}\in\mathbb{U}$, $x^a_{k+N+1|k}=f\left(x^*_{k+N|k},u^a_{k+N|k}\right)$ such that 
\begin{equation}
 \begin{aligned}
& m\left(x^a_{k+N+1|k}\right)   \leq m\left(x^*_{k+N|k}\right)-\beta_3\left(\left|x^*_{k+N|k}\right|\right) 
\le\alpha_{k+1}\\
 		\Rightarrow & x^a_{k+N+1|k}\in \Omega(\alpha_{k+1}).
 \end{aligned}	
\end{equation}
The feasible sequences at time $k+1$ can be constructed as
 \[
 \begin{aligned}
 	\bm{\mathrm u}^a_{k+1}&=\left(u^*_{k+1|k}, \cdots, u^*_{k+N-1|k}, u^a_{k+N|k}\right)\\
 	\bm{\mathrm x}^a_{k+1}&=\left(x^*_{k+1|k}, \cdots, x^*_{k+N|k}, x^a_{k+N+1|k}\right).
 \end{aligned}
 \]}

{ 
In the case of $m\left(x^*_{k+1|k}\right)< m\left(x^*_{k+N|k}\right)$, and noting   (\ref{tc}), we  have that
\begin{equation}
\begin{aligned}
	&m\left(x^*_{k+1|k}\right)< m\left(x^*_{k+N|k}\right)\le \alpha_k\\
 \Rightarrow & x^*_{k+1|k}\in \Omega(\alpha_{k}).
\end{aligned}
\end{equation}
Based on (\ref{am}), we  have that
\begin{equation}
\begin{aligned}
		\alpha_{k+1}&\le m\left(x^*_{k+1|k}\right)<m\left(x^*_{k+N|k}\right)\leq \alpha_{k}.
\end{aligned}
\end{equation}
By Assumption \textbf{(A3)}, there  exists a control $u^b_{k+1|k}\in\mathbb{U}$,
$x^b_{k+2|k}=f\left(x^*_{k+1|k},u^b_{k+1|k}\right)$ such that 
\begin{equation}\label{fk2}
 \begin{aligned}
&m\left(x^b_{k+2|k}\right)  \leq m\left(x^*_{k+1|k}\right)-\beta_3\left(\left|x^*_{k+1|k}\right|\right)\\
&~~~~~~~~~~~~~~~~~\le m\left(x^*_{k+1|k}\right)-\delta\\
&~~~~~~~~~~~~~~~~~=\alpha_{k+1}, \text{ if } \delta\leq \beta_3\left(\left|x^*_{k+1|k}\right|\right)\\
 		\Rightarrow &x^b_{k+2|k}\in \Omega(\alpha_{k+1})\subseteq \Omega(\alpha_{k}).
 \end{aligned}	
\end{equation}
Using Assumption \textbf{(A3)} again and noting (\ref{fk2}), there exists a $u^b_{k+2|k}\in\mathbb{U}$,
$x^b_{k+3|k}=f\left(x^b_{k+2|k},u^b_{k+2|k}\right)$ such that 
\begin{equation}\label{fk3}
 \begin{aligned}
&m\left(x^b_{k+3|k}\right)  \leq m\left(x^b_{k+2|k}\right)-\beta_3\left(\left|x^b_{k+2|k}\right|\right)\\
&~~~~~~~~~~~~~~~~~\le m\left(x^b_{k+2|k}\right)\le \alpha_{k+1}\\
 		\Rightarrow &x^b_{k+3|k}\in \Omega(\alpha_{k+1}).
 \end{aligned}	
\end{equation}
Combing (\ref{fk2}) and (\ref{fk3}) together, we  have that the set $\Omega(\alpha_{k+1})$ is also control invariant and hence there exists a control sequence $u^b_{k+i|k}\in \mathbb{U}$, $i\in \mathbb{I}_{0:N}$  such that 
\begin{equation}
	\begin{aligned}
	x^b_{k+1+i|k}=f\left(x^b_{k+i|k}, u^b_{k+i|k}\right)\in \Omega(\alpha_{k+1}).
\end{aligned}
\end{equation}
The feasible sequences of this case can be conducted as
 \[
 \begin{aligned}
 	\bm{\mathrm u}^b_{k+1}&=\left(u^b_{k+1|k}, \cdots, u^b_{k+N-1|k}, u^b_{k+N|k}\right)\\
 	\bm{\mathrm x}^b_{k+1}&=\left(x^*_{k+1|k}, \cdots, x^b_{k+N|k}, x^b
 	_{k+N+1|k}\right).
 \end{aligned}
 \]
In conclusion, based on the above discussion on the two cases, the proposed algorithm is recursively feasible, which completes the proof. \QEDA
 }
 
Recursive feasibility guarantees that once the proposed algorithm is feasible at time $k=0$, the optimisation algorithm associated with MPC and the closed-loop system is always feasible, which implies that all the signals in the control system (e.g., $\alpha_k$) are well defined, even as time goes to infinity.

\subsection{Asymptotic stability}
 {To establish stability in Theorem \ref{t1}, several lemmas are necessary. Lemma \ref{lem:V(x)} shows the upper bounds of the new value function $V(x,\alpha)$ in (\ref{op})  while Lemma~\ref{l3} establishes the condition for its monotonicity.}  

\begin{lem} \label{lem:V(x)}
 { Suppose that Assumption \textbf{(A3)} holds. Then, for any feasible state $x$ of the optimisation problem (\ref{op}) with $\Omega(\alpha) \subseteq \Omega$, there exists a $\mathcal{K}_\infty$ function $\beta_4(\cdot)$  such that $V(x,\alpha)\le\beta_4(|x|)$. } 
\end{lem}
{
\textit{Proof.}
We first consider a simple case, where $x\in\Omega(\alpha)$. Then,  we replace the optimisation problem (\ref{op}) by a suboptimal control problem, i.e. splitting it into $N$ one-step optimisation problems \eqref{m}. Recursively solving these $N$ one-step optimisation problems as in (\ref{m}) and using the optimal states and inputs as the feasible ones of the optimisation problem (\ref{op}), it will give that
\begin{equation}
   V(x,\alpha)\leq  N \beta_2(|x|), \forall x\in\Omega(\alpha).
\end{equation}
Therefore, $V(x,\alpha)$ is continuous at  zero as $0\in \Omega(\alpha)$. Following the similar analysis given in \cite[Chap. 2.4.2]{rawlings2017model}, we know that the set of feasible states is closed and $V(x,\alpha)$ is locally bounded by it. By applying \cite[Proposition B.25]{rawlings2017model}, the upper bound conclusion can be extended to any feasible state, which completes the proof. \QEDA}

% \begin{lem} \label{lem:V(x)}
%   Suppose that Assumption \textbf{(A4)} holds. Then, for any feasible state $x$ of the optimisation problem (\ref{op}), there exists a positive constant $\beta_4\ge\beta_1$ such that $V(x)\le\beta_4|x|^2$.  
% \end{lem}
% 	 \textit{Proof.}
% 	 First at time $k$, we replace the optimisation problem (\ref{op}) by a suboptimal control problem, i.e. splitting it into $N$ one-step optimisation problems. Recursively solving these $N$ one-step optimisation problems as in (\ref{m}). The CLF property of $m(x)$ ensures the feasibility of the optimisation problems, and, according to (\ref{mv1}),  the corresponding value functions costs, $m(x_{k+i})$, are upper bounded by 
% 	 \begin{equation} \label{eq:mk+1}
% 	     m(x_{k+1}) \le \beta_2|x_{k+i}|^2, i=0, \ldots, N-1
% 	 \end{equation}
% under the suboptimal solution. Furthermore, letting $\alpha(k):= m(x_k)$, it can be shown that $x_{k+i} \in \Omega(\alpha(k)), i=1,\ldots,N-1$. Therefore, this, together with  (\ref{eq:mk+1}), implies that there exists $\beta_4>\beta_1$ such that the corresponding cost function under the suboptimal solution is bounded by  $\beta_4|x_k|^2$. Since the optimal solution $V(x)$ must be smaller than or equal to the cost function under the suboptimal solution. So boundedness of the optimal cost function $V(x)$ in Lemma \ref{lem:V(x)} is established. 
% \QEDA

%This proof is following the spirit of those in \cite[Propositions 2.15, 2.16]{rawlings2017model}, i.e., we only need to prove that the value function $V(x_k)$ is locally bounded on $\Omega(\alpha_k)$.

\begin{lem}\label{l3}
{Suppose that  there exists a control $u_{k+N|k}\in   \mathbb{U}$ such that $x_{k+N+1|k}=f\left(x^*_{k+N|k}, u_{k+N|k}\right)\in \Omega(\alpha_{k+1})$. Then 
		\begin{align*}
			&V(x_{k+1},\alpha_{k+1})-V(x_k,\alpha_{k})\leq \ell \left(x^*_{k+N|k},u_{k+N|k},x_{k+N+1|k}\right)  -\ell \left(x^*_{k|k},u^*_{k|k},x^*_{k+1|k}\right) .
		\end{align*}}
\end{lem}
{
\textit{Proof.} We denote the optimal control and state sequences at time $k$ as 
\begin{align*}
	\bm{\mathrm u}^*_k&=\left(u^*_{k|k}, \cdots, u^*_{k+N-1|k}\right)\\
	\bm{\mathrm x}^*_k&=\left(x^*_{k|k}, \cdots, x^*_{k+N|k}\right).
\end{align*}
The value function $V(x_k,\alpha_{k})$ can be rewritten as 
\begin{equation}
	V(x_k,\alpha_{k})=\sum_{i=0}^{N-1}\ell\left(x^*_{k+i|k},u^*_{k+i|k},x^*_{k+i+1|k}\right).
	\end{equation}
If there exists a control $u_{k+N|k}\in   \mathbb{U}$ such that $x_{k+N+1|k}=f\left(x^*_{k+N|k}, u_{k+N|k}\right)\in \Omega(\alpha_{k+1})$, we could construct the feasible sequences at time $k+1$ as
\begin{align*}
	\bm{\mathrm u}_{k+1}&=\left(u^*_{k+1|k}, \cdots, u^*_{k+N-1|k},u_{k+N|k}\right)\\
	\bm{\mathrm x}_{k+1}&=\left(x^*_{k+1|k}, \cdots, x^*_{k+N|k},x_{k+N+1|k}\right).
\end{align*}
Thus, we   have that
\begin{equation}
\begin{aligned}
&V(x_{k+1},\alpha_{k+1})-V(x_{k},\alpha_{k})\\
&\leq \mathcal{J} (x_{k+1}, \bm{\mathrm u}_{k+1})-V(x_{k},\alpha_{k})\\
&=\mathcal{J} \left(x^*_{k+1|k}, \bm{\mathrm u}_{k+1}\right)-V(x_{k},\alpha_{k})\\
&=   \ell\left(x^*_{k+N|k},u_{k+N|k},x_{k+N+1|k}\right)\\
&~~~-\ell \left(x^*_{k|k},u^*_{k|k},x^*_{k+1|k}\right),
\end{aligned}
\end{equation}
which completes the proof. \QEDA}
 
{
We are now ready to prove Theorem 5 for asymptotic stability of the proposed MPC algorithm. It is established through a two stage process. In the first process, it is shown that the terminal set continuously shrinks to the origin under the constructed update rule. In the second stage, asymptotic stability is established using $V(x,0)$ as a Lyapunov function. 
}

{
\textit{Proof of  asymptotic stability in Theorem \ref{t1}.} 
Recursive feasibility established in Section~\ref{sec:RF} shows that 
\begin{equation} \label{mxk}
	\begin{aligned}
		\alpha_{k+1}\leq m\left(x^*_{k+s|k}\right)-\delta \leq m\left(x^*_{k+N|k}\right)-\delta \leq \alpha_{k}-\delta,
	\end{aligned}
\end{equation}
which implies that the sequence $\{\alpha_k\}$ is continuously decreasing until reaching the origin since $\delta_k>0$ until $x^*_{k+N|k}=0$. Consequently, once the terminal constraint set approaches the origin, one has 
\begin{equation}
%    \alpha_k=0,~\forall k\ge  {\alpha_0}/{\delta}\Rightarrow 
    x^*_{k+N|k}=0.
\end{equation}
It follows from Assumption \textbf{(A3)} that 
\begin{equation}
 \ell\left(x^*_{k+N|k},u_{k+N|k},x_{k+N+1|k}\right)=0.
 %, ~\forall k\ge  {\alpha_0}/{\delta}.   
\end{equation}
 Following Lemma \ref{l3} with $\alpha_{k+1}=\alpha_k=0$, we have 
\begin{equation}
\begin{aligned}
 & V(x_{k+1},0)-V(x_{k},0)\\
 &\leq -\ell \left(x^*_{k|k},u^*_{k|k},x^*_{k+1|k}\right)  \leq -m(x_{k|k}) \\
 &\leq -\beta_1(|x_{k|k}|).
\end{aligned}
\end{equation}
Combining Lemma \ref{lem:V(x)} and noting $V(x,\alpha)\ge m(x)\ge \beta_1(|x|)$, the closed-loop system is asymptotically stable under the proposed MPC algorithm
%after time $\alpha_0/\delta$, 
which completes the proof. \QEDA
}

{
\begin{rem}
Similar to the existing MPC theories, asymptotic stability is established also by employing an optimal value function $V(x,0)$ as a Lyapunov function candidate in our approaches. However there are three key differences in these two approaches. First, the new value function $V(x,\alpha)$ is different from the the original optimal value function obtained by optimising the original cost function (\ref{J1}) used in the existing MPC stability theories.  This is not only because the modification of the cost function by removing the cost associated with initial state, but also the contractive constraints in the online optimisation, caused by shirking $\alpha_k$. Secondly, stability of the MPC algorithm is established in the existing MPC theories by exploiting the monotonicity of the optimal value function during the whole MPC control process. However we are not able to establish the monotonicity of the value function, $V(x,\alpha)$, at the beginning of the contracive MPC algorithm proposed in this paper. Actually $V(x,\alpha)$ may increase if the terminal set shrinks quite quickly based on the proposed update rule. We rely on the contractive terminal constraints enabled by the property of OSVF being a CLF to force the system state to approach the origin. Once the terminal set shrinks to the origin, we are able to show the monotonicity of the value function $V(x,0)$. Finally, we establish stability by \emph{imposing a certain property on the modified stage cost} (i.e. OSVF to be a CLF) while the existing stability theories ensure stability by \emph{imposing a condition on the terminal cost} (i.e. the terminal cost to be a CLF).                   
\end{rem}
}
\section{Terminal cost calculation}\label{sec4}
All the assumptions used in establishing stability of MPC are quite similar to that in the existing terminal cost based MPC framework except Assumption \textbf{(A3)}. This section is devoted to {developing} a procedure for finding a terminal cost such that the corresponding OSVF $m(x)$ is a CLF. As discussed previously, the condition $m(x)$ being positive definite is opposite to the existing well-established MPC stability conditions. Furthermore, it is possible for $V_f(x)$ to be zero which includes the case of no terminal weight. Therefore, the results \cite{chen2020model} for stability analysis of MPC without terminal cost can be considered as a special case in this paper. More importantly, we are able to explore the freedom of adding an appropriate terminal cost to guarantee the stability of a MPC algorithm of concern.  

In what follows, we propose a design procedure to find a terminal cost to satisfy Assumption \textbf{(A3)} for linear systems, which can be extended to nonlinear systems by the similar approach of the current MPC framework. Linearisation of the nonlinear system at the operation condition \cite{rawlings2017model}
is one such approach,   or using another approach like linear differential inclusion  as in\cite{chen2003terminal}. It also provides more insight about the definition and the properties of OSVF $m(x)$.  
Here, we consider a linear system as 
\begin{equation}\label{ls}
	\begin{aligned}
		x^+&=Ax+Bu,\\
	\end{aligned}
\end{equation} 
with the stage cost $l(x,u)=|x|^2_Q  +|u|^2_R$. We aim to find a  terminal cost $V_f(x)=|x|^2_P$ satisfying Assumption \textbf{(A3)}. The  modified stage cost in (\ref{eq:newdefl}) can   be specified  as
\begin{equation}
	\begin{aligned}
		\ell\left(x,u,x^+\right) &=-|x|^2_P+|x|^2_Q  +|u|^2_R+|Ax+Bu|^2_P\\
		&=|x|^2_{A^TPA+Q-P}+|u|^2_{R+B^TPB}\\
		&~~~+2x^T A^TPB u\\
		&=\begin{bmatrix}
			x^T&u^T
		\end{bmatrix} M\begin{bmatrix}
			x\\u 
		\end{bmatrix},
	\end{aligned}
\end{equation}
where 	
\begin{equation} \label{eq:M}
    M\coloneqq  \begin{bmatrix}
	A^TPA+Q-P & A^TPB\\
	B^TPA &R+B^TPB
\end{bmatrix}.
\end{equation}

If $R+B^TPB>0$, solving the one-step optimisation problem (\ref{m}) gives
\begin{equation}\label{Mp}
	\begin{aligned}
		\arg\min_u \ell(x)&=-(R+B^TPB)^{-1}B^TPA x\\
		m(x)=\min_u \ell(x)&=|x|^2_{M_P},
	\end{aligned}
\end{equation}
where 
\begin{align} \label{MpD}
	M_P&\coloneqq A^TPA+Q-P-A^TPB(R+B^TPB)^{-1}B^TPA.
\end{align}
Since that $m(x)$  serves as a CLF in Assumption \textbf{(A3)}, $M_P$ should at least be positive definite. By the Schur Complement Lemma,   $R+B^TPB>0$ and $M_P>0$ is equivalent to
 $M>0$.

This result is summarised in the following lemma.

\begin{lem} \label{lem:PD}
    Consider a linear system (\ref{ls}) with a quadratic cost function. Then its OSVF $m(x)$ is positive definite if there exists a terminal cost  $V_f(x)=|x|^2_P$ with a matrix $P$ such that $M$ defined in  (\ref{eq:M}) is positive definite.
\end{lem}

By now, we only need to consider condition (\ref{mv2}). Since in this case, $m(x)$ is in a  quadratic form, condition (\ref{mv2}) can be relaxed since there exists a control $u_k$ such that $m(x_{k+1})< m(x_k)$ holds for any $x_k$.

\begin{lem}\label{l6}
$m(x_{k+1})<m(x_k)$ holds, if and only if there exist $u_k$ and $u_{k+1}$ such that $\ell(x_{k+1},u_{k+1},x_{k+2})< m(x_k)$ holds.
\end{lem} 

\textit{Proof.}
$\Leftarrow:$	If there exist $u_k$ and $u_{k+1}$ such that $\ell (x_{k+1},u_{k+1},x_{k+2})< m(x_k)$, we then have that $m(x_{k+1}) \le \ell(x_{k+1},u_{k+1},x_{k+2})< m(x_k) $.
	
$\Rightarrow:$	If  $m(x_{k+1}) < m(x_k)$, the existence of $u_k$ is obvious.  As for $u_{k+1}$, we can choose $$u_{k+1}=\arg \min_{u_{k+1}} \ell(x_{k+1}, u_{k+1}, x_{k+2}).$$ Hence,  there exist $u_k$ and $u_{k+1}$ such that $$\ell(x_{k+1},u_{k+1},x_{k+2})=m(x_{k+1})<m(x_k).$$	

This completes the proof. \QEDA

\begin{lem}
	If $R>0$,  then
	\begin{equation}
		\begin{bmatrix}
			Q-SR^{-1}S^T &A\\
			A^T&B
		\end{bmatrix}>0\Leftrightarrow \begin{bmatrix}
			Q&S &A\\
			S^T &R&0\\
			A^T&0&B
		\end{bmatrix}>0.
	\end{equation}
\end{lem}

\textit{Proof.} By using congruent transformation, we have  that
\[
\begin{aligned}
&\begin{bmatrix}
	I &-SR^{-1}&0\\
	0&I&0\\
	0&0&I
\end{bmatrix}\begin{bmatrix}
	Q &S&A\\
	S^T&R&0\\
	A^T&0&B
\end{bmatrix} \begin{bmatrix}
	I &0&0\\
	-R^{-1}S^T&I&0\\
	0&0&I
\end{bmatrix}\\
&~~~~~~~~=\begin{bmatrix}
	Q-SR^{-1}S^T &0&A\\
	0&R&0\\
	A^T&0&B
\end{bmatrix}.	
\end{aligned}
 \] 
This completes the proof. \QEDA

\begin{thm}
 If there exists a terminal cost $P$ and two control gains $K_1$ and $K_2$ such that the following bilinear matrix inequality (BMI) is satisfied
\begin{equation}\label{bmi}
\begin{bmatrix}
		M &\begin{bmatrix}
			(A+BK_1)^T&  K^T_2\\
			0&0
		\end{bmatrix}  M \\
		M\begin{bmatrix}
			(A+BK_1)&0\\K_2 &0  
		\end{bmatrix} &  M
	\end{bmatrix} >0
\end{equation}
then $m(x)$ is a CLF for system (\ref{ls}).
\end{thm}	

\textit{Proof.} Note that if condition (\ref{bmi}) is satisfied,  it follows that \(M>0\), which further guarantees  $R+B^TPB>0$ and $M_P>0$ by using the Schur Complement Lemma. $R+B^TPB>0$ implies that the one-step optimisation problem (\ref{m}) is well-defined.    

Letting   $u_k= K_1x_k$ and  $u_{k+1}= \tilde K_2x_{k+1}$,  
a quick calculation gives
\begin{equation}
	\begin{aligned}
	&\ell(x_{k+1},u_{k+1},x_{k+2}) \\
	&=\begin{bmatrix}
			x_{k}^T(A+BK_1)^T&x^T_k K^T_2
		\end{bmatrix}M\begin{bmatrix}
			(A+BK_1)x_{k}\\K_2x_k 
		\end{bmatrix}	\\
		&=x_{k}^T\begin{bmatrix}
			(A+BK_1)^T&  K^T_2
		\end{bmatrix}M\begin{bmatrix}
			A+BK_1\\K_2  
		\end{bmatrix}x_{k}	
	\end{aligned}
\end{equation}
where $K_2=\tilde K_2(A+BK_1)$. Using  the  Schur Complement Lemma and Lemma \ref{l6}, noting $R+B^TPB>0$, the condition $\ell(x_{k+1},u_{k+1},x_{k+2})< m(x_k)$   is satisfied if the following matrix inequality holds

\[
	\begin{aligned}
		& \begin{bmatrix}
			(A+BK_1)^T&  K^T_2
		\end{bmatrix}M\begin{bmatrix}
			A+BK_1\\K_2  
		\end{bmatrix}< M_P \\
		\Leftrightarrow &  \begin{bmatrix}
			(A+BK_1)^T&  K^T_2
		\end{bmatrix}  M   M^{-1}  M\begin{bmatrix}
			A+BK_1\\K_2  
		\end{bmatrix}< M_P \\
		\Leftrightarrow & \begin{bmatrix}
			M_P &\begin{bmatrix}
				(A+BK_1)^T&  K^T_2
			\end{bmatrix}  M \\
			  M\begin{bmatrix}
				A+BK_1\\K_2  
			\end{bmatrix} &  M
		\end{bmatrix}> 0\\
		\Leftrightarrow & ~(\ref{bmi}).
	\end{aligned}
\]

This completes the proof. \QEDA

%
%\begin{equation}
%\begin{aligned}
%&A_K\coloneqq A+BK, ~M_P\coloneqq A^TPA+Q-P-A^TPB(R+B^TPB)^{-1}B^TPA>0\\
%	&A_K^TM_P A_K<M_P\\
%	\Leftrightarrow & M^{-1}_P A_K^TM_P A_KM^{-1}_P<M^{-1}_P\\
%	\Rightarrow &\begin{bmatrix}
%		M_P^{-1} & M^{-1}_P A_K^T \\
%		A_KM^{-1}_P &M_P^{-1}
%		\end{bmatrix}>0
%	\end{aligned}
%\end{equation}

\begin{rem}
The proposed BMI (\ref{bmi}) can be solved by using the PENLAB package \cite{fiala2013penlab}, which is an open source software package implemented in MATLAB for matrix inequalities.	
\end{rem}

\section{Illustrative examples}\label{sec5}

\subsection{A first order example}
This section is to illustrate the significance of the proposed new stability analysis approach in this paper and highlight the main differences with the existing stability theory of MPC by using a simple first order system
\begin{equation}
	\begin{aligned}
		x^+=ax+bu,~b\neq0
	\end{aligned}
\end{equation}
with $l(x,u)=qx^2+ru^2$ and $V_f(x)=px^2$, where the lowercase letters  represent scalars.

More specifically, this example is used to highlight the following two points. 
\begin{enumerate}
    \item The proposed MPC algorithm may be stable with a zero terminal cost or a negative terminal cost with the new condition.
    \item The new stability condition is complementary to the existing terminal cost based MPC stability conditions. 
\end{enumerate}
   
It follows from the definition of the augmented stage cost that $\ell(x,u,x^+)=-px^2+qx^2+ru^2+p(ax+bu)^2$ for this example. The corresponding OSVF $m(x)$ is ready to be calculated.  For this first order linear system, according to Lemma~\ref{lem:PD}, the positive definite condition in Assumption \textbf{(A3)}    reduce to
\begin{equation}\label{fo1}
	\begin{aligned}
	 r+b^2p&>0\\
     a^2p+q-p-\frac{a^2p^2b^2}{r+b^2p}&>0.
	\end{aligned}
\end{equation}

Furthermore, Lemma~\ref{l6} shows that the CLF condition (\ref{mv2}) is met for any positive definite OSVF that satisfies condition (\ref{fo1}). 

Inspecting condition (\ref{fo1})  gives
\begin{equation}\label{lowbq}
	\begin{aligned}
q>(1-a^2)p+\frac{a^2p^2b^2}{r+b^2p} =\frac{1}{b^2}\left(z_1+\frac{a^2r^2}{z_1}-\left(1+a^2\right)r\right),
	\end{aligned}
\end{equation}
where $z_1\coloneqq r+b^2p>0$.  Eq.(\ref{lowbq}) gives a lower bound of $q$ by
\begin{equation}\label{qmin}
\begin{aligned}
		q> \underline q &\coloneqq\frac{2|ar|-(1+a^2)r}{b^2}\\
		&=\left\{\begin{aligned}
		&\frac{-(|a|-1)^2r}{b^2},\text{ if }r>0\\
		&0,~~~~~~~~~~~~~~~~~\text{ if }r=0\\
		&\frac{-(|a|+1)^2r}{b^2},\text{ if }r<0.
	\end{aligned}\right.
\end{aligned}
\end{equation}  
It shall be highlighted that Eq.(\ref{qmin}) reveals  $q$ is possible to be negative or zero if  $r\neq0$,.

Next, for the purpose of comparison, we present the conventional MPC stability conditions for this simple MPC problem. First, all the weightings are required to be at least positive semi-definite, i.e., $q>  0$, $r\ge   0$, and $p\ge  0$. Secondly, the current condition states that the conventional MPC for this unconstrained system is stable if the famous Fake Algebraic Riccati Equation (FARE) \cite{BitGevWer90} holds
\begin{equation}\label{lowbq2}
	\begin{aligned}
		&b^2 p^2+(1-b^2q-a^2)p-q\ge 0\\
		&\Leftrightarrow q\le \frac{1}{b^2}\left(z_2+\frac{a^2r^2}{z_2}-\left(1+a^2\right)r\right)\\
	\end{aligned}
\end{equation}
where $z_2=r+b^2p\ge r$.

If we focus on the  right-side functions of (\ref{lowbq}) and (\ref{lowbq2}), when $p$ is large enough,  both functions go to a linear one $q=p-\frac{ra^2}{b^2}$, which is parallel with but below the line $q=p$ if $r\neq 0$.

In what follows, we will make an explicit comparison of the stability conditions and regions between the proposed new approach and the existing terminal weight based MPC framework, in terms of the control weight $r>0$, $r<0$, and $r=0$.

\subsubsection{Case with control weight $r>0$}
\begin{figure}[ht]
	\footnotesize
	\begin{center} 
		\begin{tabular}{c}
			\includegraphics[width=0.35\textwidth]{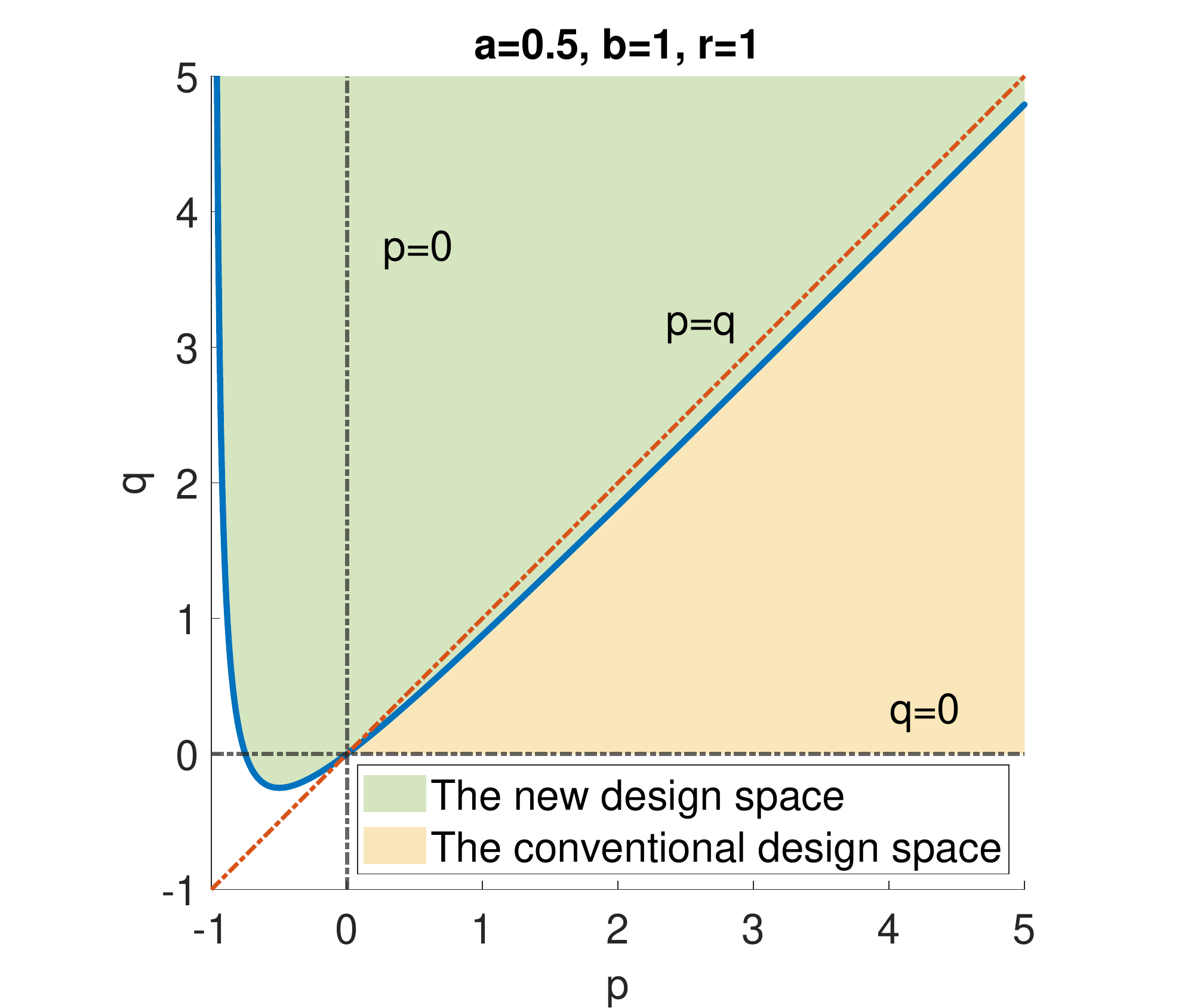} \\
			a)\\
			\includegraphics[width=0.35\textwidth]{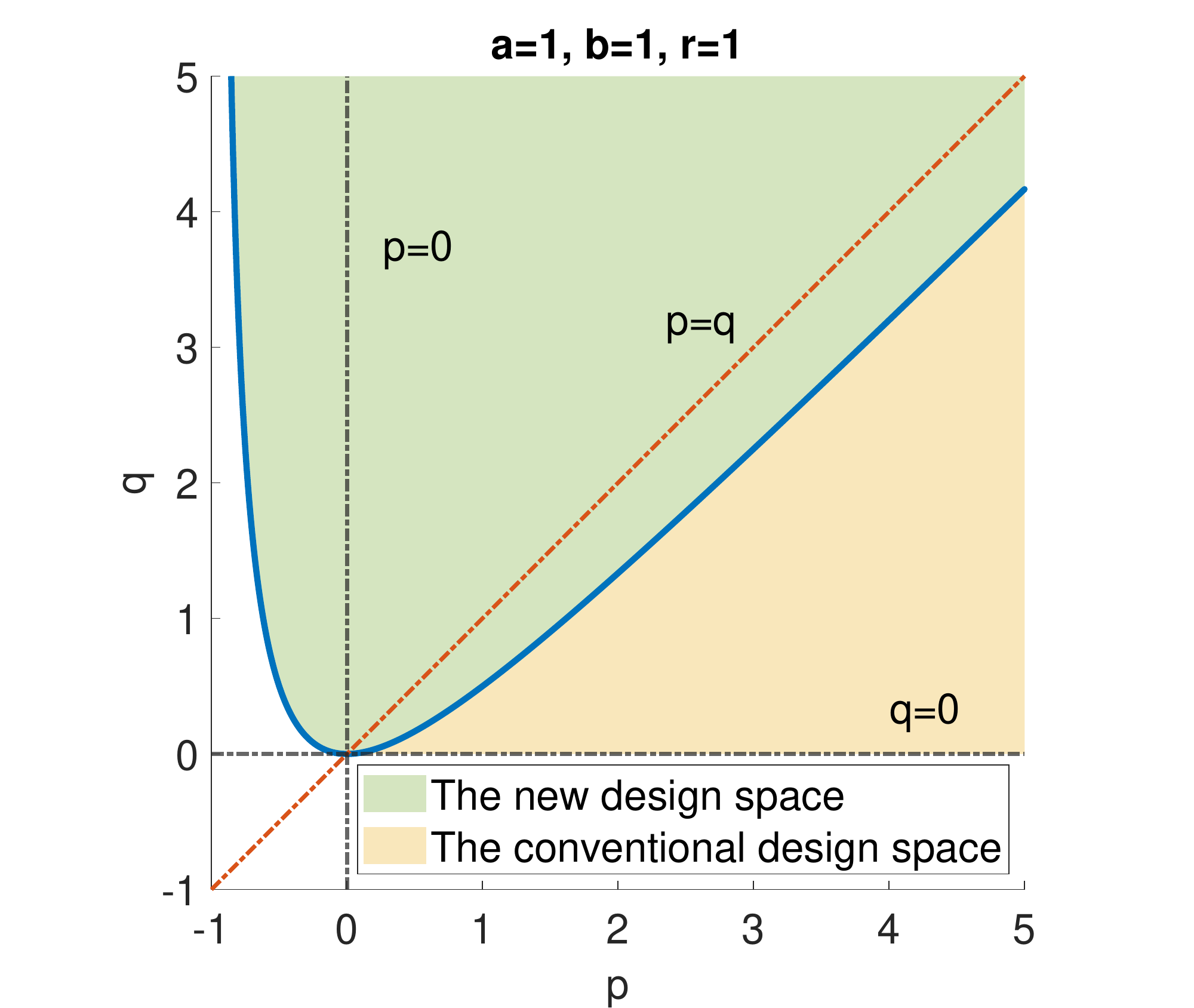} \\
			b)\\
		  \includegraphics[width=0.35\textwidth]{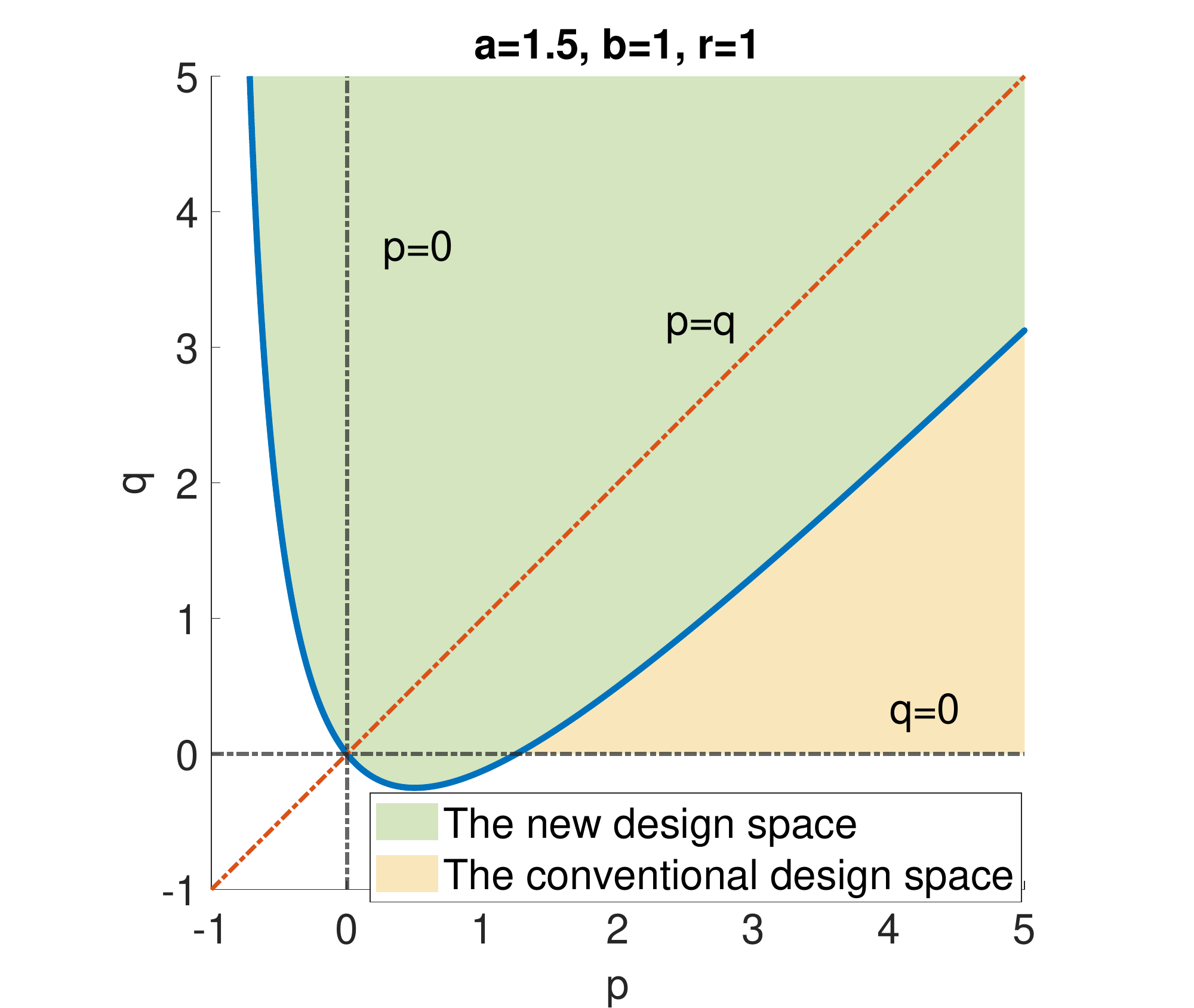}  \\
		 c)
		\end{tabular}
		\caption{Design spaces of the proposed  and conventional MPC  for different open-loop systems when $r>0$. a) Stable, i.e. $|a|<1$. b) Marginally stable, i.e. $|a|=1$. c) Unstable, i.e. $|a|>1$. } \label{sr1}
	\end{center}
\end{figure}

We only consider the case that  $r=1$, as the same results hold for a different $r$ by normalizing state and terminal weighting with $q/r$ and $p/r$. Regarding stability condition (\ref{qmin}), there are three cases to discuss, $|a|<1$, $|a|=1$, and $|a|>1$, which correspond to open-loop stable, marginally stable, and unstable systems, respectively. The stability regions defined by the new condition (\ref{qmin}) are given by Fig. \ref{sr1}. These regions actually define the possible parameter combinations of the system and cost functions including the stage and terminal cost so we refer them as the feasible \emph{design spaces}.  It is interesting to see  when $|a|<1$, the proposed MPC could be stable even when both $p$ and $q$ are negative, which covers the situation widely studied in EMPC. 

It shall be  noted that the blue curve corresponds to the solution of the Riccati equations under different stage costs, i.e. $q$. Therefore, to ensure stability, the existing MPC stability condition requires that the terminal cost $p$ is larger than or equal to the solution of the Riccati equation. This is shown by the region on the right side of the blue curve; that is, the system under MPC is stable when the corresponding terminal cost $p$ covers the cost-to-go. In contrast, the stability region given by our new stability condition lies on the left side of the blue curve so covers the design space where the terminal weight $p$ is less than the Riccati solution under the corresponding stage and control weights $q$ and $r$.          

It is observed that there is no overlap between the stability regions derived by the existing terminal weight based MPC framework and the new condition proposed in this paper. Furthermore, by combining the stability regions yielded by these two methods, we are able to significantly extend the MPC design space (i.e. in terms of possible terminal costs) where stability could be guaranteed. Therefore, our approach is complementary to the existing MPC framework.

Finally, is is also observed that both the lines $p=0$ and $p=q$ are covered by the proposed MPC stability condition as shown in Fig \ref{sr1}. $p=0$ and $p=q$ are both the special cases of MPC where there is no terminal cost, but with slightly different Lyapunov functions. To well illustrate the difference, we take an example when $N=2$.
\begin{itemize}
    \item If $p=0$, the optimisation in (\ref{op}) reduces to  
    \begin{equation}
	\begin{aligned}
		&\min_{u_{k|k},~u_{k+1|k}}	  \left(qx^2_{k|k}+ru^2_{k|k}\right) +\left(qx^2_{k+1|k}+ru^2_{k+1|k}\right)
		\\\text{ s.t. } &  x_{k|k}=x_k\\
		&x_{k+i+1|k}=ax_{k+i|k}+bu_{k+i|k},~i=0,1\\
		& x_{k+2|k} \in \Omega(\alpha_k),  
	\end{aligned}
\end{equation}
where  $M_P=a^2p+q-p-a^2p^2b^2(r+b^2p)^{-1}=q.$
 \item If $p=q$, the optimisation in (\ref{op}) becomes 
 \begin{equation}
	\begin{aligned}
		&\min_{u_{k|k},~u_{k+1|k}}	  \left(qx^2_{k|k}+ru^2_{k|k}\right)\\
		&~~~~~~~~~~~~~~~~~+\left(qx^2_{k+1|k}+ru^2_{k+1|k}\right)\\
		&~~~~~~~~~~~~~~~~~+px^2_{k+2|k}\\
		=&\min_{u_{k|k},~u_{k+1|k},~u_{k+2|k}}	 \sum_{i=0}^{2} \left(qx^2_{k+i|k}+ru^2_{k+i|k}\right)
		\\\text{ s.t. } &  x_{k|k}=x_k\\
	&x_{k+i+1|k}=ax_{k+i|k}+bu_{k+i|k},~i=0,1\\
		& x_{k+2|k} \in \Omega(\alpha_k),  
	\end{aligned}
\end{equation}
where  $M_P=a^2p+q-p-a^2p^2b^2(r+b^2p)^{-1}=a^2q-a^2q^2b^2(r+b^2q)^{-1}= a^2qr(r+b^2q)^{-1}$.  The equivalent cost function holds as the  decision variable $u_{k+2|k}$ would never affect the control and state sequences. 
\end{itemize}
Therefore, they correspond to the different choices of the augmented stage cost and the associated OSVF. This also implies different terminal constraints are used in the algorithm and different Lyapunov functions employed in establishing stability. This special case also provides more insight into the proposed new approach.

%\begin{figure*}[ht]
%	\footnotesize
%	\begin{center} 
%		\begin{tabular}{cc}
%			\includegraphics[width=0.4\textwidth]{a05r1} &
%			\includegraphics[width=0.4\textwidth]{a1r1} \\
%			a)&b)\\
%			\multicolumn{2}{c}{\includegraphics[width=0.4\textwidth]{a15r1}} \\
%			\multicolumn{2}{c}{c)} 
%		\end{tabular}
%		\caption{Design spaces of the proposed  and conventional MPC  for different open-loop systems when $r>0$. a) Stable, i.e. $|a|<1$. b) Marginally stable, i.e. $|a|=1$. c) Unstable, i.e. $|a|>1$. } \label{sr1}
%	\end{center}
%\end{figure*}

\subsubsection{Case with control weight $r=0$}
\begin{figure}[ht]
	\footnotesize
	\begin{center} 
		\begin{tabular}{c}
			\includegraphics[width=0.35\textwidth]{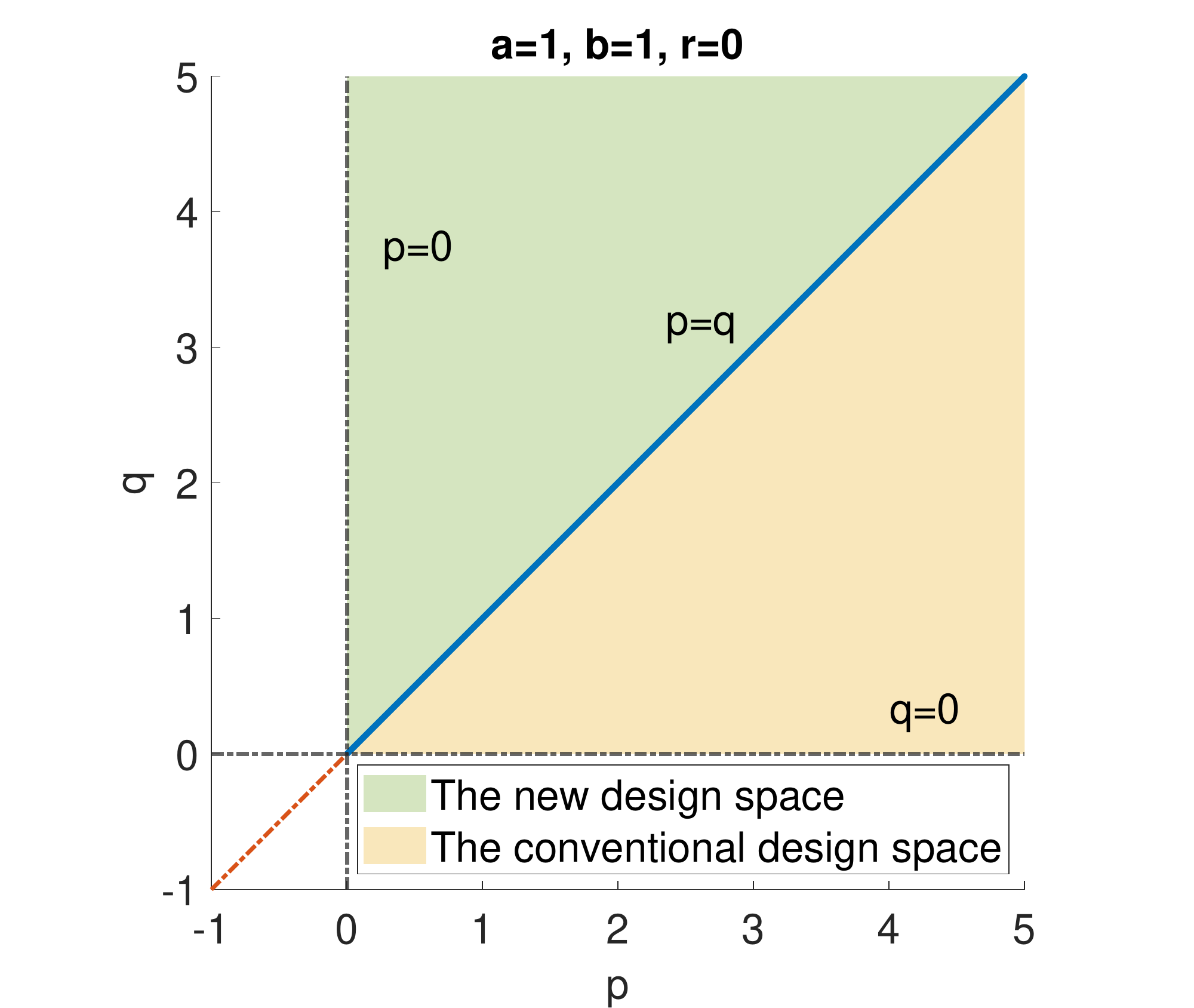} \\
			a)\\
			 \includegraphics[width=0.35\textwidth]{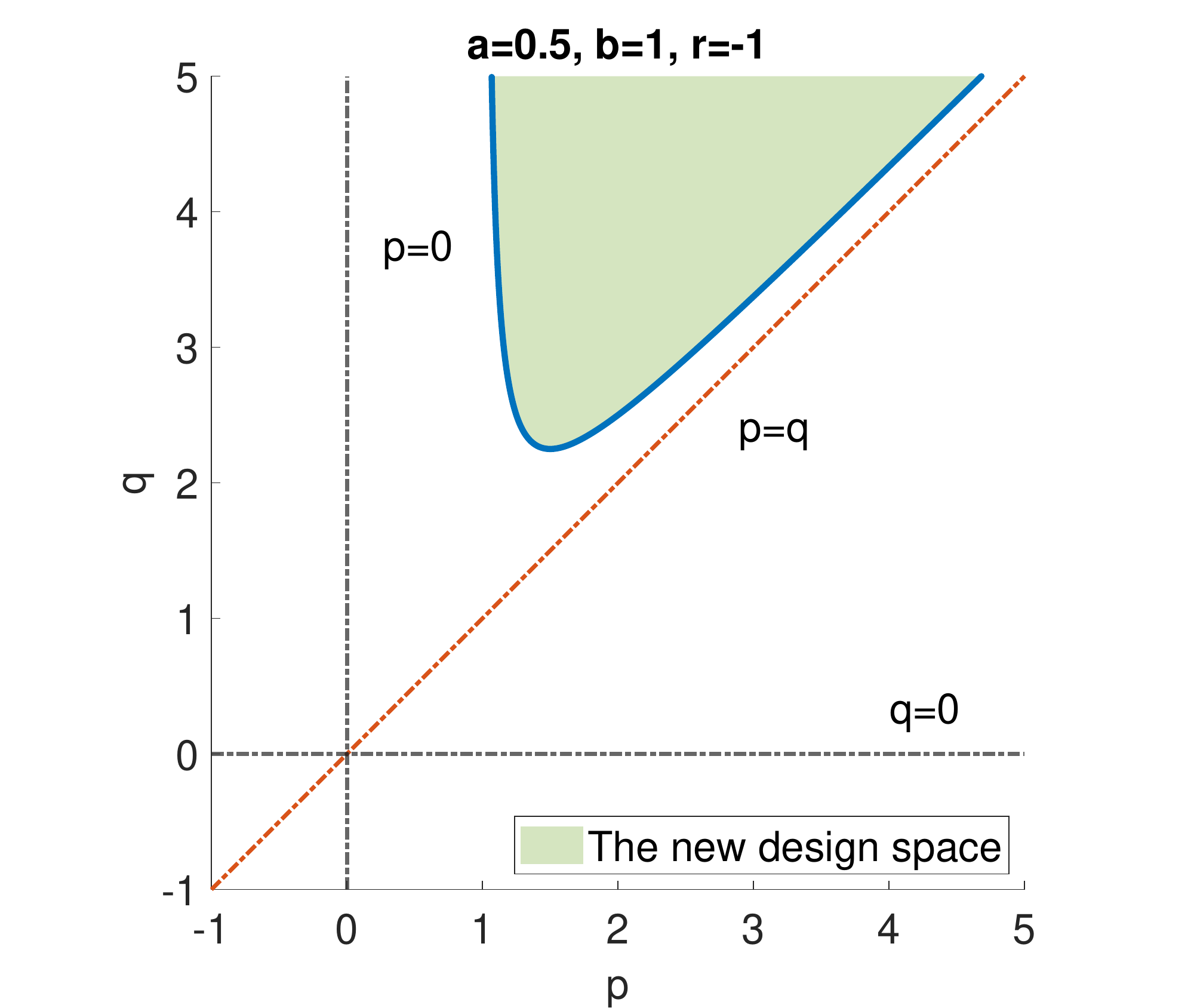} \\
			b)
		\end{tabular} 
		\caption{Design spaces of the proposed  and conventional MPC.   a) $r=0$. b) $r<0$.} \label{sr2}
	\end{center}
\end{figure}

The case of $r=0$ is also interesting as both methods occupy the perfect half of the whole design space, as shown by Fig. \ref{sr2} a). This is because when $r=0$, condition (\ref{lowbq}) is equivalent to $q>p>0$ while  condition (\ref{lowbq2}) is $0 < q\le  p$.

%\begin{figure*}[ht]
%	\footnotesize
%	\begin{center} 
%		\begin{tabular}{cc}
%			\includegraphics[width=0.4\textwidth]{a1r0} & \includegraphics[width=0.4\textwidth]{a05r_1} \\
%			a)&b)\\
%		\end{tabular} 
%		\caption{Design spaces of the proposed  and conventional MPC.   a) $r=0$. b) $r<0$.} \label{sr2}
%	\end{center}
%\end{figure*}

\subsubsection{Case with control weight $r<0$}

In this case, in a similar fashion as $r>0$, we fix $r=-1$ and the result is given by Fig. \ref{sr2} b). It is worth noting that the design space of the conventional MPC framework is empty since the control weight is normally required to be positive definite.

In summary, for the cases of the control weight $r<0$ or $r=0$, we are able to draw the same conclusion from Fig \ref{sr2} about the proposed new stability condition, and its comparison and connection with the existing terminal weight based MPC framework, as in the case of $r>0$. More simply, for each pair of state and control weights $q,r$, the existing MPC stability theory ascertains that the corresponding MPC algorithm with the terminal cost $p$ in the yellow region is stable while our stability condition states that the corresponding MPC algorithm with the terminal cost $p$ in the light green area is also stable. Thus, we are able to state that an MPC algorithm is stable with any terminal cost from the combined areas of these two stability approaches. This substantially relaxes the stability requirement on MPC algorithms and increase the freedom in designing and tuning the key parameters involved in MPC. This numerical example and detailed comparison give more insight into the proposed new augmented stage cost method and its links with the existing MPC stability theory. 

\subsection{Simple nonlinear system study}
\begin{figure}[ht]
	\footnotesize
	\begin{center} 
\includegraphics[width=0.35\textwidth]{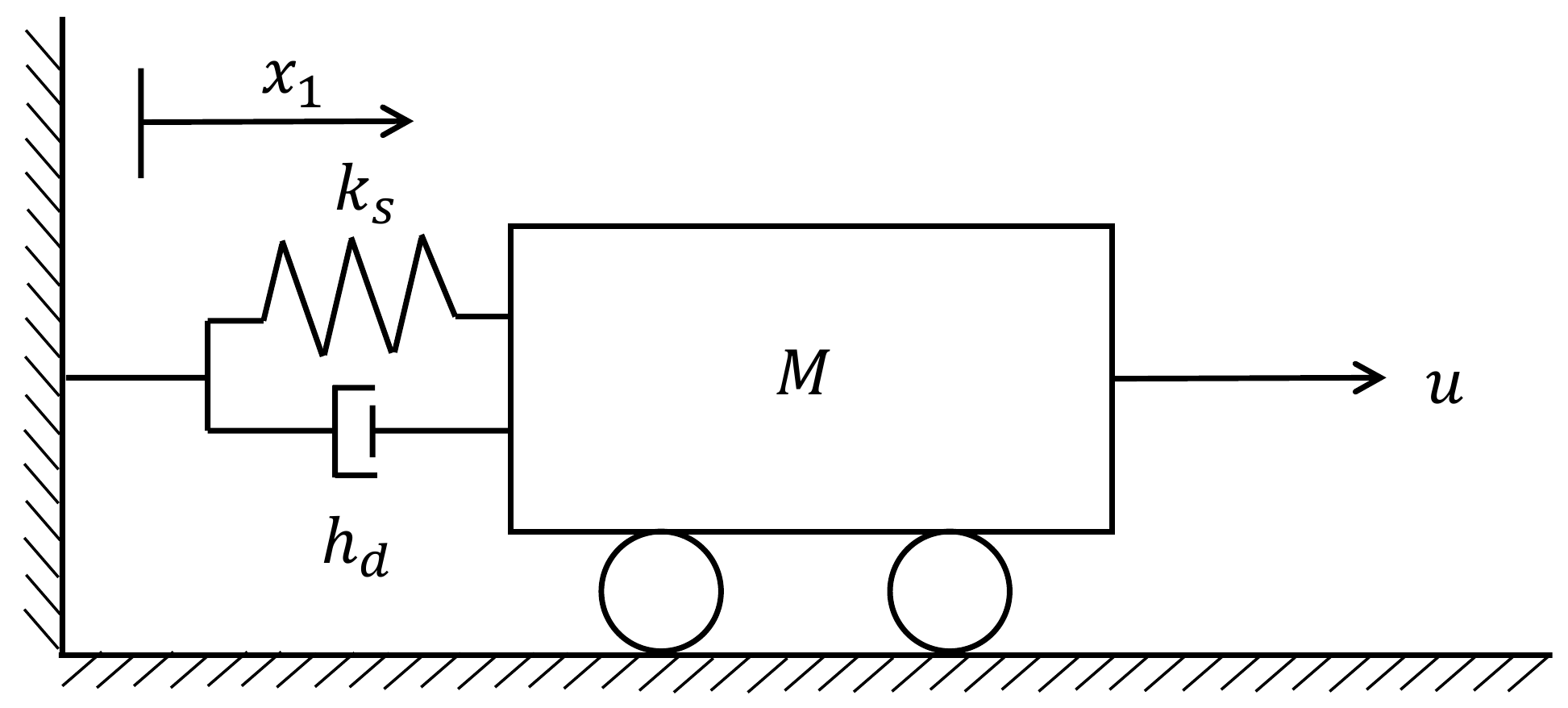}  		
\caption{Cart and spring example} \label{car}
	\end{center}
\end{figure}

This section is to apply the new proposed MPC algorithm with a contractive terminal constraint on a nonlinear system and compare it with the existing terminal weight based MPC approach.  It aims to highlight their differences in terminal regions and  total running costs. Here, these two MPC laws are applied to a benchmark system \cite{magni2003robust}, i.e,  a cart with a mass moving on a plane, as shown in Fig. \ref{car}. The nonlinearity in this example comes from the elasticity in the spring, $k_s=k_0e^{-x_1}$, where $x_1$ is the  displacement of the carriage from the equilibrium position associated with the external force $u=0$. 
%Besides, a damper with damping factor $h_d$ affects the system in a resistive way. For easy understanding, 
The system is described by
\begin{equation}
	\begin{aligned}
		\dot {x}_1&=x_2\\
		\dot {x}_2&=-\frac{k_0}{M}e^{-x_1}x_1-\frac{h_d}{M}x_2+\frac{u}{M},
	\end{aligned}
\end{equation}
where $x_2$ is its velocity and $h_d$ is the damping coefficient. The system is discretised into a discrete-time system using the Euler approximation with sampling time of 0.4sec, given by  
 \begin{equation}
 	\begin{aligned}
 		x_{1,k+1}&=x_{1,k}+0.4x_{2,k}\\
 		x_{2,k+1}&=-0.132e^{-x_{1,k}}+0.56x_{2,k}+0.4u_k,
 	\end{aligned}
 \end{equation}
where the specific system parameters are all given in  \cite{magni2003robust}. The state and input constrains are $\mathbb{X}=[-2,2]\times[-3,3]$ and  $\mathbb{U}=[-4,4]$ and the initial state is $x_0=[-2,1]^T$.  The linearization at the origin gives the controllable pair 
\[
A=\begin{bmatrix}
	1 & 0.4\\
	-0.132&0.56
\end{bmatrix},~B=\begin{bmatrix}
	0\\
	0.4
\end{bmatrix}.
\]
The prediction horizon is chosen to be short enough, $N=3$, due to the fast calculation time by the motion control system.  The stage and terminal costs are both given in a quadratic form, i.e., $l(x,u)=|x|^2_Q  +|u|^2_R$ and $V_f(x)=|x|^2_P$, where
\[Q=\begin{bmatrix}
	2 & 0\\
	0&4
\end{bmatrix},~R=1.\]

\begin{figure}[ht]
	\footnotesize
	\begin{center} 
			\includegraphics[width=0.35\textwidth]{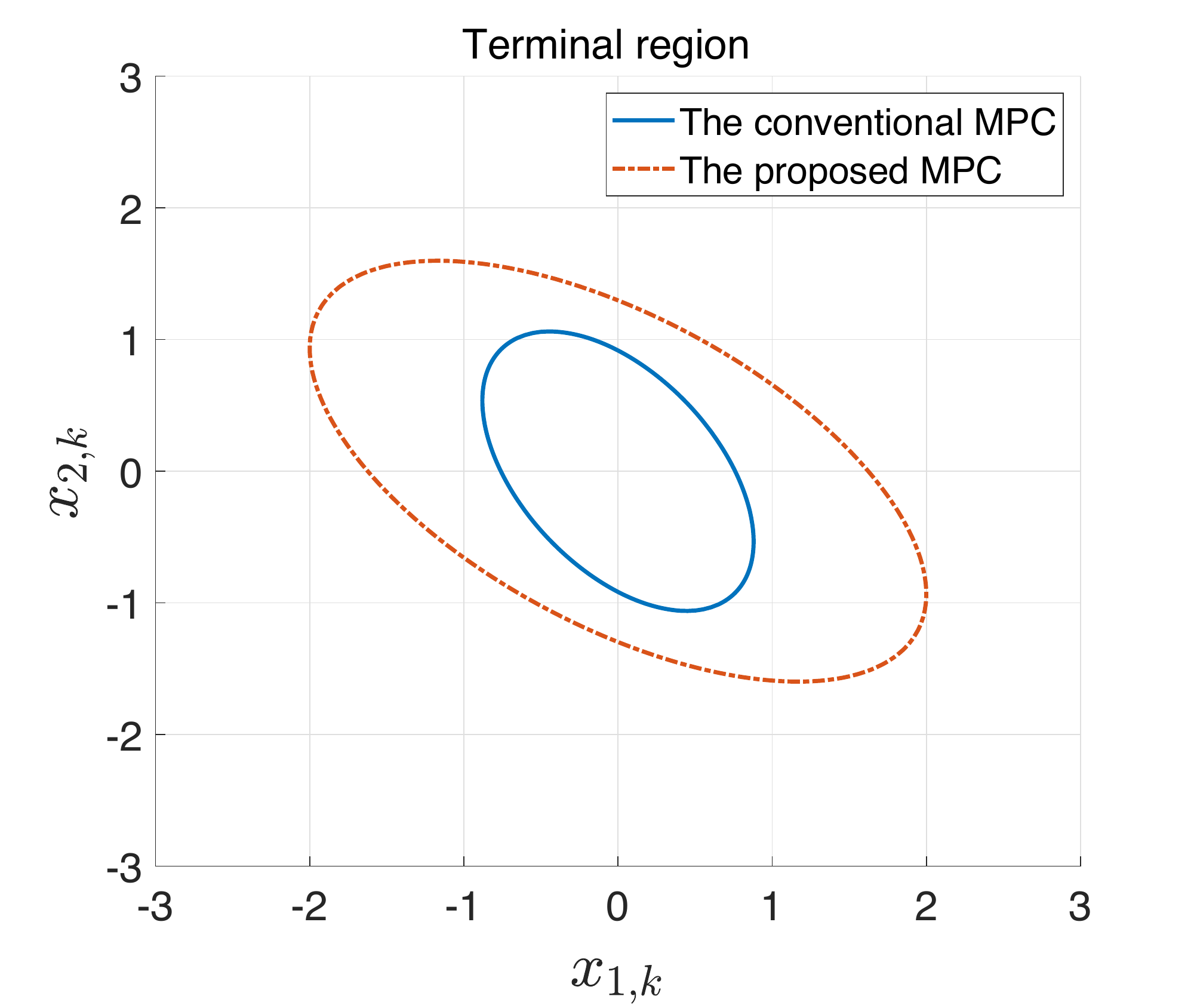} \\
		\caption{Terminal regions of the proposed  and conventional MPC.}\label{simtr}
	\end{center}
\end{figure}

\begin{figure}[ht]
	\footnotesize
	\begin{center} 
		\begin{tabular}{c}
			\includegraphics[width=0.35\textwidth]{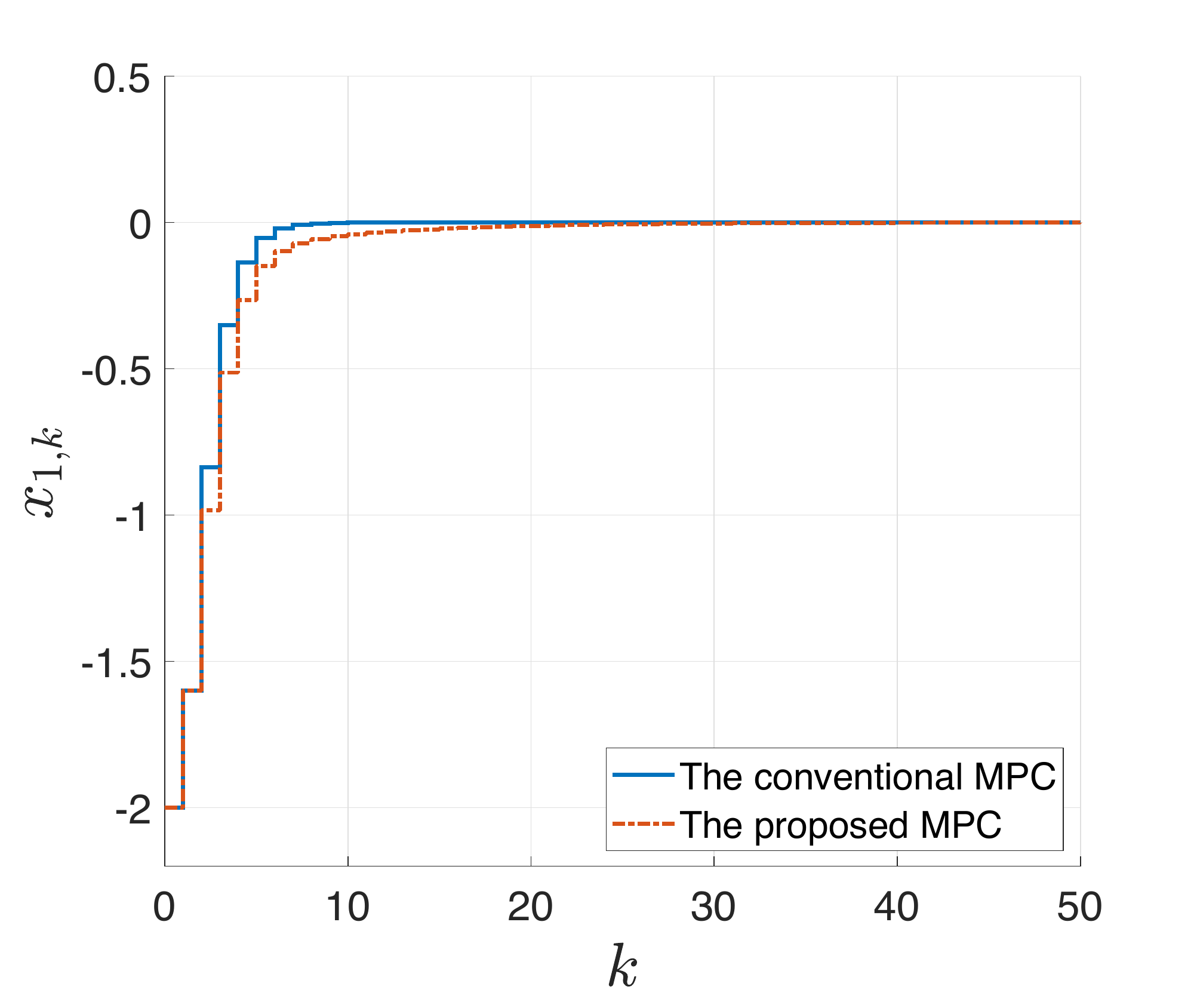} \\
			a)\\
			\includegraphics[width=0.35\textwidth]{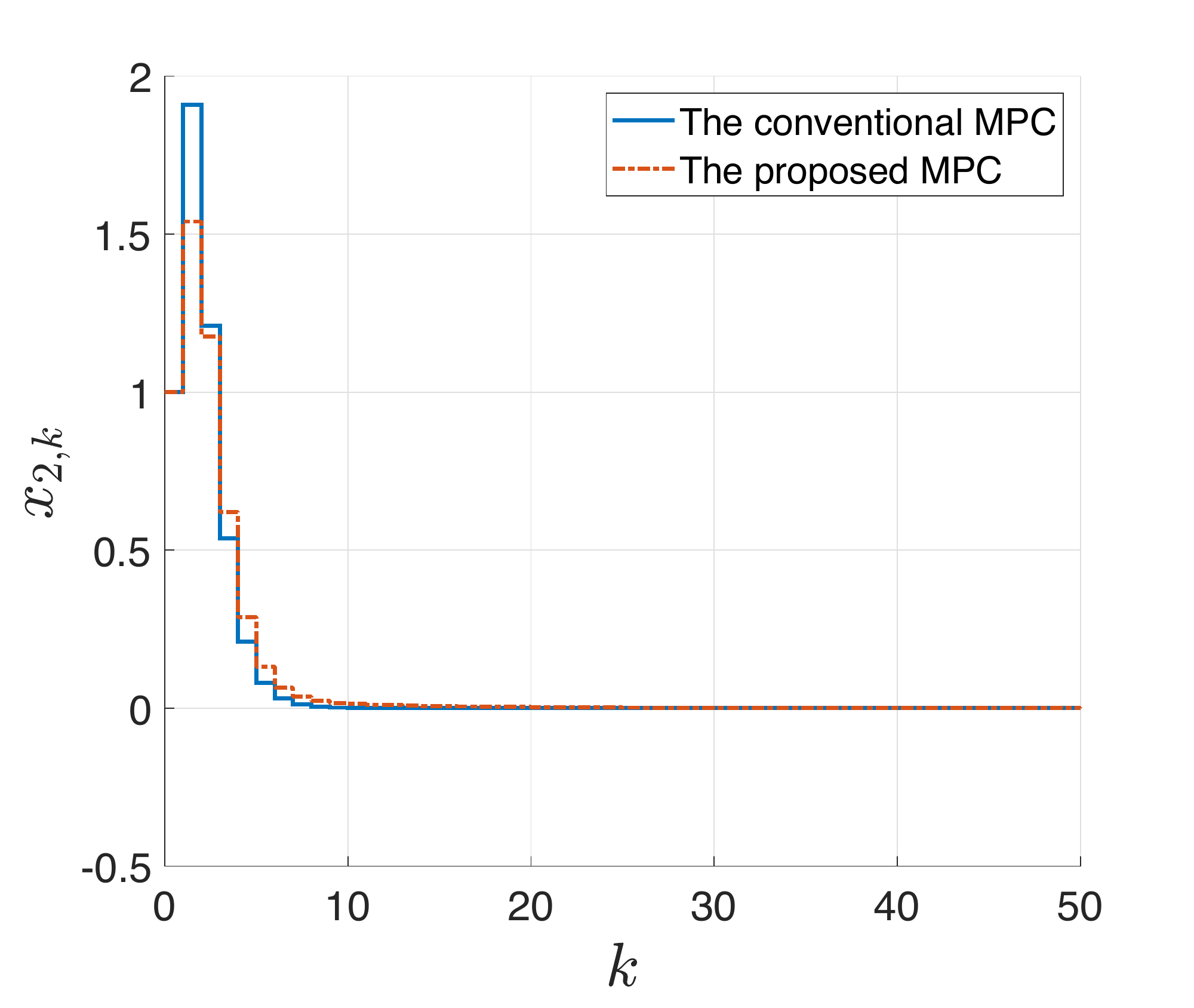} \\
			b)\\
			\includegraphics[width=0.35\textwidth]{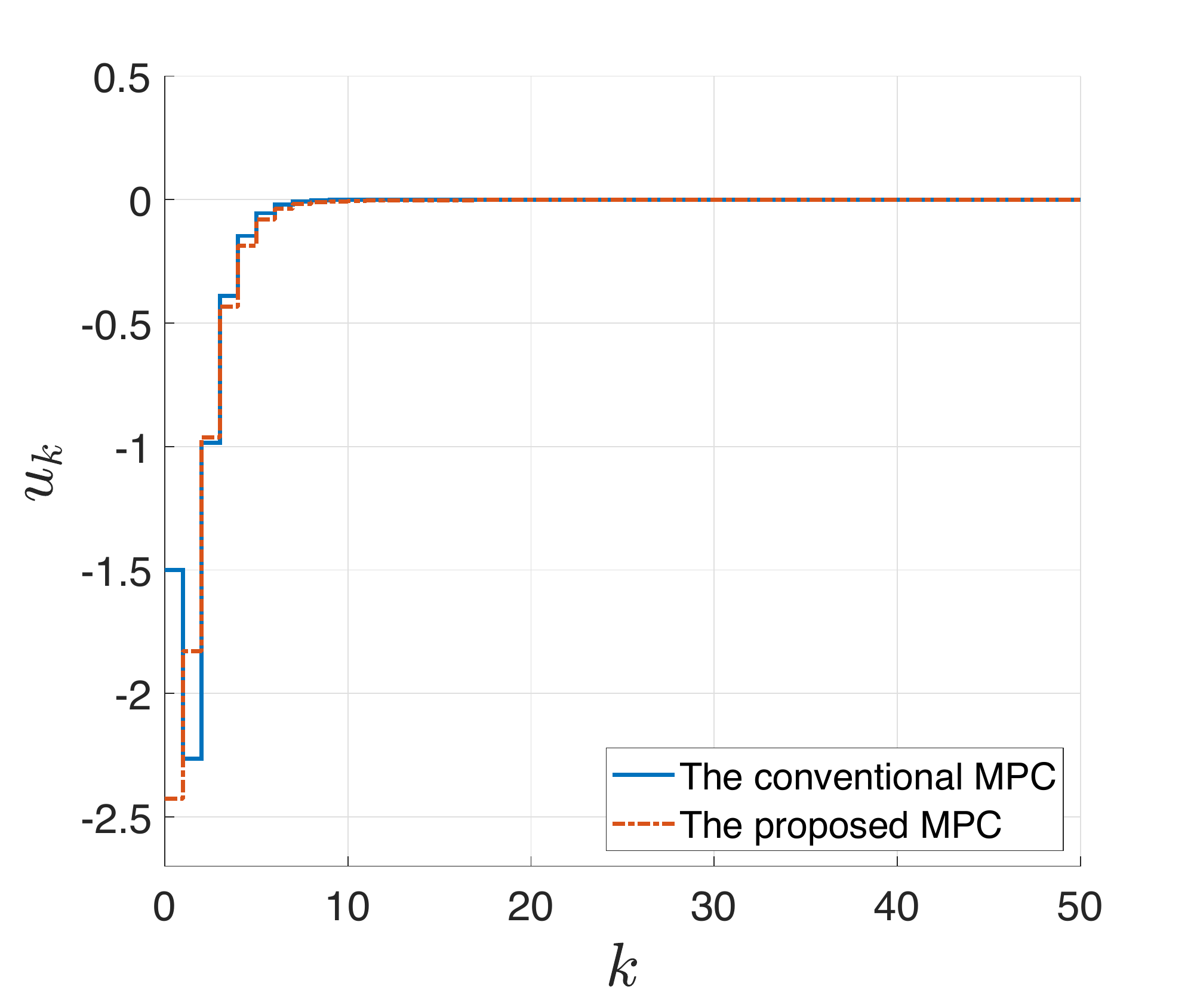}\\  
			c)
		\end{tabular}
		\caption{States and inputs of the proposed  and conventional MPC. a) $x_{1,k}$. b) $x_{2,k}$. c) $u_k$. }\label{simsi} 
	\end{center}
\end{figure}

In what follows, we will calculate the terminal costs and regions and assess their stability and performance using these two methods respectively. 
\begin{itemize}
	\item  We first consider the conventional MPC design. The terminal cost is obtained by the Riccati equation for the linearization
	\[P=\begin{bmatrix}
		 10.9153  &  4.5604\\
    4.5604  &  7.5023
	\end{bmatrix}\] 
	and the terminal region is given by
	\[
	\Omega=\left\{x\in R^2: x^TPx\le \alpha\right\}
	\] 
	where $\alpha=6.3076$ is calculated by the optimization method introduced in \cite{chen2003terminal}. 
	\item  As for the proposed MPC design, following the procedure in Section \ref{sec4}, the terminal cost is obtained by solving the BMI  (\ref{bmi})
	\[P=\begin{bmatrix}
    3.5249 &  -0.3522\\
   -0.3522 &   1.5731
	\end{bmatrix}\] 
    and the corresponding terminal region is given by
    \[
    \Omega=\left\{x\in R^2: x^TM_Px\le \alpha\right\}
    \] 
    where $M_P$ is computed by (\ref{MpD}) and $\alpha=5.4823$ is also calculated by  the same method in \cite{chen2003terminal}. 
\end{itemize}
The terminal regions are given by Fig. \ref{simtr}, where the region of the proposed MPC covers that of the conventional one.  The states and inputs are given by Fig. \ref{simsi}. It shall be highlighted that different from stability regions plotted in the parameter space as in Figs. \ref{sr1} and \ref{sr2}, the terminal regions are given in state space in Fig. \ref{simsi}.  It is worth noting that the terminal cost $P$ calculated by our approach  does not satisfy the condition (\ref{eq:original}), i.e.  covering the cost to go as required in the existing terminal weight based MPC framework.   

Although the original stage costs used in these two MPC algorithms are the same, the final cost functions are different since their terminal costs are different. The total running cost is used for a fair comparison of the performance of these two methods, that is,  \[J_{\text{run}}=\sum_{k=0}^{T_s} l\left(x_k,u_k\right)=\sum_{k=0}^{T_s} \left(|x_k|^2_Q  +|u_k|^2_R\right),\]
where $T_s$ is sufficiently large such that the system reaches its steady state, which is chosen as 50 sec in this study. The total running cost of the proposed MPC algorithm is 47.2148, which is less than that of the conventional one at 49.1587.

Therefore, for this example, comparing with the existing terminal weight based MPC scheme, the proposed MPC algorithm has a larger terminal region and exhibits slightly better performance.  

%\begin{figure}[ht]
%	\footnotesize
%	\begin{center} 
%		\begin{tabular}{cc}
%			\includegraphics[width=0.4\textwidth]{x1} &
%			\includegraphics[width=0.4\textwidth]{x2} \\
%			a)&b)\\
%			\multicolumn{2}{c}{\includegraphics[width=0.4\textwidth]{u}} \\
%			\multicolumn{2}{c}{c)} 
%		\end{tabular}
%		\caption{States and inputs of the proposed  and conventional MPC. a) $x_{1,k}$. b) $x_{2,k}$. c) $u_k$. }\label{simsi} 
%	\end{center}
%\end{figure}

\section{Conclusion}\label{sec6}
Stability is an important property of any control method since it helps to ensure the safety of a system. Comparing with other optimisation based methods, one distinctive feature of finite horizon optimisation based control is  that it is possible to establish its stability under certain conditions. However, despite its huge success, the current stability conditions are still conservative, which restricts the design space of MPC algorithms. This paper explores the design space where the current terminal weight based MPC framework is not applicable. We propose a completely complementary stability condition in this unexplored design space. We show that, by constructing a contractive terminal constraint, rather than a fixed terminal constraint as in the current MPC framework, we are able to guarantee stability in a significant region that has not been covered previously. To achieve this, we construct an augmented stage cost using a terminal cost. That is, the augmented stage cost consists of a nominal stage cost and a rotated terminal cost. It is shown that feasibility and the stability of the proposed MPC algorithm can be established under the condition that the one-step optimal value function of the augmented stage cost is a CLF. Due to the nature of the proposed stability condition, it complements to the existing stability theory and is not intended to replace the conventional stability conditions. The stability region of MPC can be significantly enlarged by combining the new condition with the existing ones; e.g. guarantee stability under a much wider range of combination of state, control and terminal weights. The work in this paper only concerns a basic MPC regulation problem without disturbance or uncertainty.   

Our future work includes extending the proposed approach to other areas including robust MPC, stochastic MPC or economic MPC where stability has been established with the current terminal weight based MPC approach and comparing their differences.   

% e condition that augmented stage cost is       

% As a stability guaranteed optimisation method,  Model Predictive Control (MPC) has been continuously showing its significance after its born, especially in the era of optimisation domination. However, enormous successful MPC applications which do not rigorously follow the conventional stability conditions, e.g., the MPC without any terminal cost, imply the huge gap between the conventional MPC theory and practice.

%  In this paper, we have further bridged such a gap by  establishing a complementary condition to guarantee the stability of MPC, although the proposed new algorithm slightly modifies the MPC setup by introducing a time-varying contractive terminal set. It has been proved that the proposed algorithm is recursively feasible and asymptotically stable, as also provided by the conventional MPC. On the other hands, as a complementary approach, the proposed algorithm would never replace the conventional MPC, but significantly enlarge the design space of MPC variants.

\section*{Acknowledgments}
This work was supported by the UK Engineering and Physical Sciences Research Council (EPSRC) Established Career Fellowship ``Goal-Oriented Control Systems: Disturbance, Uncertainty and Constarints'' under the grant number EP/T005734/1.

%Bibliography
\bibliographystyle{unsrt}

\end{document}